\documentclass[11pt]{article}

\usepackage{amsmath,amssymb,amsfonts,amsthm}
\usepackage{graphicx}
\usepackage{array}
\usepackage{url}
\usepackage{xcolor}
\usepackage[colorlinks=true,allcolors=blue]{hyperref}

\newtheorem{theorem}{Theorem}[section]
\newtheorem{proposition}[theorem]{Proposition}
\newtheorem{lemma}[theorem]{Lemma}
\newtheorem{corollary}[theorem]{Corollary}

\newtheorem{remark}[theorem]{Remark}

\newcommand{\mult}{\mathrm{mult}}
\newcommand{\ch}{\mathrm{ch}}
\newcommand{\Vir}{\mathrm{Vir}}

\newcommand{\black}[1]{\textcolor{black}{#1}}

\title{\bfseries
Branching rules for winding subalgebras of the affine Kac--Moody algebras
$A^{(1)}_1$ and $A^{(2)}_2$
}

\author{Duc--Khanh Nguyen}

\date{}

\begin{document}

\maketitle

\begin{abstract}
We study branching problems for affine Kac--Moody algebras. Unlike the
finite-dimensional case, an affine Kac--Moody algebra may contain proper
subalgebras isomorphic to itself, such as winding subalgebras obtained by
rescaling the loop parameter. We investigate the restriction of integrable
highest-weight representations to such subalgebras. The restriction remains
integrable and decomposes into irreducible components with finite multiplicities, encoded by pairs of highest weights. We show that this set is closed under addition, extending a result of Brion and Knop to the affine setting. We also give a partial description of this set and provide explicit results for types $A^{(1)}_1$ and $A^{(2)}_2$.
\end{abstract}




\section{Introduction}\label{Sect1}
A central problem in representation theory is to understand how an irreducible module over a Lie algebra \(\mathfrak{g}\) decomposes when restricted to a Lie subalgebra \(\dot{\mathfrak{g}}\).
When \(\mathfrak{g}\) and \(\dot{\mathfrak{g}}\) are finite-dimensional semisimple Lie algebras, irreducible modules are parametrized by dominant integral weights \(P_+\) and \(\dot{P}_+\), respectively.
For \(\Lambda \in P_+\), the restriction of the irreducible \(\mathfrak{g}\)-module \(L(\Lambda)\) to \(\dot{\mathfrak{g}}\) decomposes as
\begin{equation}\label{eq:1}
L(\Lambda)
=
\bigoplus_{\lambda \in \dot{P}_+}
\dot{L}(\lambda)^{\mult_{\Lambda,\dot{\mathfrak{g}}}(\lambda)},
\end{equation}
and the determination of the multiplicities \(\mult_{\Lambda,\dot{\mathfrak{g}}}(\lambda)\) is known as the branching problem.

This framework encompasses many classical situations.
Diagonal embeddings give rise to tensor product multiplicities, including the Littlewood--Richardson coefficients for \(\mathfrak{gl}_n(\mathbb{C})\), while restriction to a Cartan subalgebra $\mathfrak{h}$ recovers the Kostka coefficients.
Beyond individual multiplicities, the structure of the support
\[
\Gamma(\mathfrak{g},\dot{\mathfrak{g}})
=
\bigl\{(\Lambda,\lambda)\in P_+\times\dot{P}_+
\mid \mult_{\Lambda,\dot{\mathfrak{g}}}(\lambda)\neq 0\bigr\}
\]
is of independent interest: it is a finitely generated semigroup whose associated cone is polyhedral \cite{Brion2}.
In the diagonal case, this cone is the Horn cone, which has been extensively studied; see, for instance,
\cite{Brion, Kumarsv, Brion2, Ressayreinv, BelkaleKumar}.

In this paper, we investigate analogous branching phenomena for affine Kac--Moody algebras.
Let \(\mathfrak{g}\) be an affine Kac--Moody algebra and \(L(\Lambda)\) an integrable highest weight \(\mathfrak{g}\)-module.
Finite-multiplicity decompositions of the form~\eqref{eq:1} occur when \(\dot{\mathfrak{g}}\) is a Cartan subalgebra, is diagonally embedded, or is a winding subalgebra in the sense of Kac--Wakimoto~\cite{KacWaki}.
While the first two cases have been widely studied, especially in connection with tensor product decompositions of affine and symmetrizable Kac--Moody algebras
\cite{Ressayretensor, RessayreKumar, KumarBrown, Peter},
our focus here is on the winding case.
This setting is closely related to tensor product multiplicities and to the proof of Frenkel's conjecture \cite{KacWaki, Frenkel1981}, and is distinguished by the fact that the winding subalgebras \(\mathfrak{g}[u]\) are isomorphic to \(\mathfrak{g}\) but embedded in a highly nontrivial way.

Let \(P_+\) and \(\dot{P}_+\) denote the sets of dominant integral weights of \(\mathfrak{g}\) and \(\mathfrak{g}[u]\), respectively.
By \cite{KacWaki}, for \(\Lambda \in P_+\) the restriction of \(L(\Lambda)\) to \(\mathfrak{g}[u]\) remains integrable and decomposes as
\[
L(\Lambda)
=
\bigoplus_{\lambda \in \dot{P}_+}
\dot{L}(\lambda)^{\mult_{\Lambda,\mathfrak{g}[u]}(\lambda)}.
\]
We denote by
\[
\Gamma(\mathfrak{g},\mathfrak{g}[u])
=
\bigl\{(\Lambda,\lambda)\in P_+ \times \dot{P}_+
\mid \mult_{\Lambda,\mathfrak{g}[u]}(\lambda)\neq 0\bigr\}
\]
the support of this decomposition. Our first main result is the following theorem.

\begin{theorem}[Theorem \ref{M1}]\label{m1}
The set $\Gamma(\mathfrak{g},\mathfrak{g}[u])$ is a subsemigroup of $\mathfrak{h}^* \times \mathfrak{h}^*$.
\end{theorem}

Let \(\delta\) be the basic imaginary root of \(\mathfrak{g}\).
Describing \(\Gamma(\mathfrak{g},\mathfrak{g}[u])\) is equivalent to determining, for each
\((\Lambda,\lambda)\in P_+\times\dot{P}_+\), the set
\[
\mathcal{B}(\Lambda,\lambda)
=
\{\, b\in\mathbb{C}\mid \dot{L}(\lambda+b\delta)\subset L(\Lambda)\,\}.
\]
Let \(P_u(\Lambda)\) denote the set of \(\lambda\) such that
\((\Lambda,\lambda)\in\Gamma(\mathfrak{g},\mathfrak{g}[u])\).
Then \(\mathcal{B}(\Lambda,\lambda)\) is nonempty if and only if
\(\lambda\in P_u(\Lambda)+\mathbb{C}\delta\).
In this case, we define \(b_{\Lambda,\lambda,u}\) as the maximal element of
\(\mathcal{B}(\Lambda,\lambda)\) and define \(h_{\Lambda,\lambda}^{[u]}\) as in~\eqref{hAa}.
The structure of \(\mathcal{B}(\Lambda,\lambda)\) is described in the next theorem.
\begin{theorem}[Theorem \ref{BAa}]\label{BAa.intro}
Let $\mathfrak{g}$ be an affine Kac--Moody algebra of type $X_N^{(r)}$. Fix $u \in \mathbb{Z}_{>1}$ such that $u \equiv 1 \pmod r$. Let $\Lambda \in P_+$ and $\lambda \in P_u(\Lambda)+\mathbb{C}\delta$. Then:
\begin{enumerate}
    \item[1.] $b_{\Lambda,\lambda,u}-(\mathbb{Z}_{\geq 0}\setminus \{1\}) \subset \mathcal{B}(\Lambda,\lambda) \subset b_{\Lambda,\lambda,u}-\mathbb{Z}_{\geq 0}$.
    \item[2.] If, in addition, $h_{\Lambda,\lambda}^{[u]} \neq 0$, then $\mathcal{B}(\Lambda,\lambda) = b_{\Lambda,\lambda,u}-\mathbb{Z}_{\geq 0}$.
\end{enumerate}
\end{theorem}

Since \(\mathfrak{h}\subset \mathfrak{g}[u]\), any pair
\((\Lambda,\lambda)\in\Gamma(\mathfrak{g},\mathfrak{g}[u])\) satisfies that
\(\lambda\) is a weight of \(L(\Lambda)\).
We denote the corresponding weight space by
\[
L(\Lambda)_\lambda
=
\{\, v\in L(\Lambda)\mid hv=\lambda(h)v \text{ for all } h\in\mathfrak{h}\,\},
\]
and define
\[
\Gamma(\mathfrak{g},\mathfrak{h})
=
\{\,(\Lambda,\lambda)\in P_+\times\mathfrak{h}^*
\mid L(\Lambda)_\lambda\neq 0\,\}.
\]

Let \(P(\Lambda)\) be the set of weights of \(L(\Lambda)\).
For any \(\lambda\in P(\Lambda)+\mathbb{C}\delta\), the set of \(b\in\mathbb{C}\) such that
\(\lambda+b\delta\) is a weight of \(L(\Lambda)\) is of the form
\(b_{\Lambda,\lambda}-\mathbb{Z}_{\geq 0}\), for a uniquely determined
\(b_{\Lambda,\lambda}\).
Thus \(\Gamma(\mathfrak{g},\mathfrak{h})\) is completely determined by these numbers.
For \(\mathfrak{g}\) of types \(A^{(1)}_1\) and \(A^{(2)}_2\), they are explicitly computed in
Propositions~\ref{maxLamdaandnk1} and~\ref{maxLamdaandnk2}.

To formulate our results for \(\Gamma(\mathfrak{g},\mathfrak{g}[u])\) in these cases,
let \(\Lambda_0\) denote the fundamental weight and \(\alpha=\alpha_1\) the simple root.
Up to tensoring by a one-dimensional \(\mathfrak{g}\)-module, any integrable highest weight
\(\mathfrak{g}\)-module is of the form \(L(\Lambda)\) with
\[
\Lambda=m\Lambda_0+\frac{j\alpha}{2},
\qquad m,j\in\mathbb{Z}_{\geq 0},
\]
where \(j\le m\) in type \(A^{(1)}_1\) and \(j\le \tfrac{m}{2}\) in type \(A^{(2)}_2\).

Theorems~\ref{b=b.a11.intro} and~\ref{b=b.a22.intro} below extend \cite[Theorem~2.2]{KacWaki} for \(A^{(1)}_1\) and \(A^{(2)}_2\) from level one to arbitrary positive level, respectively.

\begin{theorem}[Theorem \ref{b=b.a11}]\label{b=b.a11.intro}
Let $\mathfrak{g}$ be the affine Kac--Moody algebra of type $A^{(1)}_1$. Let $\Lambda = m\Lambda_0 + \frac{j\alpha}{2} \in P_+$ and let $\lambda = m'\Lambda_0 + \frac{j'\alpha}{2} \in \dot{P}_+$.
\begin{enumerate}
    \item[1.] If there exists $b\in \mathbb{C}$ such that $(\Lambda,\lambda+b\delta)$ belongs to $\Gamma(\mathfrak{g},\mathfrak{g}[u])$, then $j'-j \in 2\mathbb{Z}$ and $m'=m$.
    \item[2.] If moreover
    \begin{enumerate}
        \item[a.] $j \leq j' \leq um-j$ and $u$ is even; or
        \item[b.] $j \leq j' \leq um - (m-j)$ and $u$ is odd,
    \end{enumerate}
    then $b_{\Lambda,\lambda,u} = b_{\Lambda,\lambda}$.
\end{enumerate}
\end{theorem}

\begin{theorem}[Theorem \ref{b=b.a22}]\label{b=b.a22.intro}
Let $\mathfrak{g}$ be the affine Kac--Moody algebra of type $A^{(2)}_2$. Let $\Lambda = m\Lambda_0 + \frac{j\alpha}{2} \in P_+$ and let $\lambda = m'\Lambda_0 + \frac{j'\alpha}{2} \in \dot{P}_+$.
\begin{enumerate}
    \item[1.] If there exists $b\in \mathbb{C}$ such that $(\Lambda,\lambda+b\delta)$ belongs to $\Gamma(\mathfrak{g},\mathfrak{g}[u])$, then $m'=m$.
    \item[2.] If moreover
    \begin{enumerate}
        \item[a.] $j \leq j'$; and
        \item[b.] $j' \in \frac{m(u-1)}{2}-j + \bigl(2\mathbb{Z}_{\geq 0} \cup \mathbb{Z}_{<0}\bigr)$,
    \end{enumerate}
    then $b_{\Lambda,\lambda,u} = b_{\Lambda,\lambda}$.
\end{enumerate}
\end{theorem}

Let $\mathcal{A}_u(\Lambda)$ be the set of all $\lambda$ for which Theorems~\ref{b=b.a11.intro} and~\ref{b=b.a22.intro} apply; see~\eqref{AuA.11}, \eqref{AuA.22}. The saturated set of the support $\Gamma(\mathfrak{g},\mathfrak{g}[u])$ is defined by
\begin{equation*}
\tilde{\Gamma}(\mathfrak{g},\mathfrak{g}[u])
=
\bigl\{(\Lambda,\lambda) \in P_+ \times \dot{P}_+ \mid \lambda \in \Lambda+Q,\ \dot{L}(N\lambda) \subset L(N\Lambda) \text{ for some integer } N>1 \bigr\}.
\end{equation*}

The following result is a corollary of Theorems~\ref{BAa.intro}, \ref{b=b.a11.intro}, and~\ref{b=b.a22.intro}.

\begin{corollary}[Corollary \ref{satured}]\label{satured.intro}
Let $\mathfrak{g}$ be an affine Kac--Moody algebra of type $A^{(1)}_1$ or $A^{(2)}_2$. Fix $u \in \mathbb{Z}_{>1}$ ($u$ is an odd integer in the case of $A^{(2)}_2$). Let $\Lambda \in P_+$ and let $\lambda \in \mathcal{A}_u(\Lambda) \cap (\Lambda+Q)$. For all $b \in \mathbb{C}$, we have
\begin{enumerate}
    \item[1.] $(\Lambda,\lambda+b\delta) \in \tilde{\Gamma}(\mathfrak{g},\mathfrak{g}[u]) \Leftrightarrow d(\Lambda,\lambda+b\delta) \in \Gamma(\mathfrak{g},\mathfrak{g}[u])$ for all $d \in \mathbb{Z}_{\geq 2}$.
    \item[2.] If, in addition, ${h}_{\Lambda,\lambda}^{[u]} \neq 0$, then
    \[
    (\Lambda,\lambda+b\delta) \in \tilde{\Gamma}(\mathfrak{g},\mathfrak{g}[u]) \Leftrightarrow (\Lambda,\lambda+b\delta) \in \Gamma(\mathfrak{g},\mathfrak{g}[u]).
    \]
\end{enumerate}
\end{corollary}

This article is organized as follows.
In Section~\ref{Sect2}, we review basic facts about affine Kac--Moody algebras.
Section~\ref{Sect3} is devoted to branching with respect to Cartan subalgebras; in particular, we prove that the set $\Gamma(\mathfrak{g},\mathfrak{h})$ forms a semigroup.
In Section~\ref{Sect4}, we introduce winding subalgebras and study the associated branching problem for general affine Kac--Moody algebras. The main results of this section are Theorems~\ref{m1} and~\ref{BAa.intro}.
Section~\ref{A1A2} treats the explicit cases of types $A^{(1)}_1$ and $A^{(2)}_2$. The main results there are Theorems~\ref{b=b.a11.intro}, \ref{b=b.a22.intro}, and Corollary~\ref{satured.intro}.
Finally, in Appendix~\ref{app}, we present the coset construction for winding subalgebras in both the untwisted and twisted affine Kac--Moody settings.


\section{Preliminaries}\label{Sect2}
In this section, we recall basic results on affine Kac--Moody algebras from \cite{Carter, Kac}.

\subsection{Affine Cartan matrices}

Let $I=\{0,\dots,l\}$. An \black{affine Cartan matrix} is a matrix
$A=(a_{ij})_{i,j\in I}$ which is indecomposable of corank~$1$, satisfies
$a_{ii}=2$, $-a_{ij}\in\mathbb{Z}_{\ge0}$ for $i\neq j$, $a_{ij}=0$
if and only if $a_{ji}=0$, and $Au=0$ for some column vector $u$
with positive integer entries. Let
$a={}^t(a_0,\dots,a_l)$ and $c=(c_0,\dots,c_l)$ be relatively prime
positive integer vectors such that
\[
Aa={}^t(0,\dots,0), \qquad cA=(0,\dots,0).
\]
The \black{Coxeter number} and the \black{dual Coxeter number} of $A$ are
\[
\mathrm{h}=\sum_{i\in I} a_i, \qquad \mathrm{h}^{\vee}=\sum_{i\in I} c_i.
\]
\subsection{Realization of a generalized Cartan matrix}

Let $(\mathfrak{h},\Pi,\Pi^{\vee})$ be a \black{realization} of $A$, where
$\mathfrak{h}$ is a $\mathbb{C}$-vector space of dimension $l+2$,
$\Pi^{\vee}=\{h_0,\dots,h_l\}\subset\mathfrak{h}$ and
$\Pi=\{\alpha_0,\dots,\alpha_l\}\subset\mathfrak{h}^*$ are linearly independent,
and $\alpha_i(h_j)=a_{ji}$.
The element $K=\sum_{i\in I} c_i h_i$ is the \black{canonical central element},
and $\delta=\sum_{i\in I} a_i \alpha_i$ is the \black{basic imaginary root}.
The \black{scaling element} $d\in\mathfrak{h}$ satisfies
$\alpha_0(d)=1$ and $\alpha_i(d)=0$ for $i>0$.
The \black{fundamental weights} $\Lambda_i$ ($i\in I$) are defined by
$\Lambda_i(h_j)=\delta_{ij}$ and $\Lambda_i(d)=0$.
Set $\rho=\sum_{i\in I}\Lambda_i$.
Then $\{h_0,\dots,h_l,d\}$ and $\{\alpha_0,\dots,\alpha_l,\Lambda_0\}$
are bases of $\mathfrak{h}$ and $\mathfrak{h}^*$, respectively.

\subsection{Affine Kac--Moody algebras}

Let $\mathfrak{g}=\mathfrak{g}(A)$ be the \black{affine Kac--Moody algebra}
associated with $A$. Then $\mathfrak{h}$ is a \black{Cartan subalgebra},
and $\Pi$ and $\Pi^{\vee}$ are the \black{sets of simple roots} and
\black{simple coroots}, respectively. We have the triangular decomposition
\[
\mathfrak{g}=\mathfrak{n}_-\oplus\mathfrak{h}\oplus\mathfrak{n}_+,
\]
where $\mathfrak{n}_-$ and $\mathfrak{n}_+$ are the \black{negative} and
\black{positive subalgebras}.
An affine Kac--Moody algebra has \black{type} $X_N^{(r)}$, $r=1,2,3$
\cite{Kac}. In particular, the types $A_1^{(1)}$ and $A_2^{(2)}$ correspond
to the generalized Cartan matrices
\[
\begin{pmatrix}
2 & -2\\
-2 & 2
\end{pmatrix},
\qquad
\begin{pmatrix}
2 & -4\\
-1 & 2
\end{pmatrix},
\]
respectively. The transpose ${}^t\!A$ is again an affine Cartan matrix;
the type of $\mathfrak{g}({}^t\!A)$ is denoted by $X_{N}^{(r^\vee)}$.
\subsection{Weyl group}

Let $Q=\mathbb{Z}\Pi$ be the \black{root lattice} and
$Q_+=\mathbb{Z}_{\ge0}\Pi$. Define a \black{partial order} on $\mathfrak{h}^*$
by $\lambda\ge\mu$ if $\lambda-\mu\in Q_+$. For $i\in I$, set
$\alpha_i^\vee=\frac{a_i}{c_i}\alpha_i$. Define a sublattice $M\subset Q$ by
\[
M=
\begin{cases}
\displaystyle \bigoplus_{i=1}^l \mathbb{Z}\alpha_i, & r^\vee=1,\\[0.4em]
\displaystyle \bigoplus_{i=1}^l \mathbb{Z}\alpha_i^\vee, & r^\vee>1.
\end{cases}
\]

Let $(\cdot|\cdot)$ be the \black{standard invariant bilinear form} on
$\mathfrak{h}^*$, given by
\[
(\alpha_i|\alpha_j)=\frac{c_i}{a_i}a_{ij}, \quad
(\alpha_i|\Lambda_0)=\frac{1}{a_0}\alpha_i(d), \quad
(\Lambda_0|\Lambda_0)=0,
\]
for all $i,j\in I$.

For $\alpha\in\mathfrak{h}^*$ define $t_\alpha\in GL(\mathfrak{h}^*)$ by
\[
t_\alpha(\lambda)
=\lambda+\lambda(K)\alpha-
\left(\lambda+\frac{\lambda(K)\alpha}{2}\,\middle|\,\alpha\right)\delta.
\]

Let $W$ be the \black{Weyl group} of $\mathfrak{g}$, generated by the
\black{fundamental reflections} $s_i$ $(i\in I)$, where
\[
s_i(\lambda)=\lambda-\lambda(h_i)\alpha_i.
\]
Then
\begin{equation}\label{Weylgroup}
W \cong t_M \rtimes \overline{W},
\end{equation}
with $t_M=\{t_\alpha\mid\alpha\in M\}$ and
$\overline{W}=\langle s_1,\dots,s_l\rangle$.

\subsection{Realization of affine Kac--Moody algebras}

Let $\overline{\mathfrak{g}}$ be a simple Lie algebra with Lie bracket
$[\,\cdot\,,\,\cdot\,]_0$ and normalized invariant form $(\cdot|\cdot)_0$
\cite{Carter}. Define
\begin{equation*}
\hat{\overline{\mathfrak{g}}}
=\mathbb{C}[t,t^{-1}]\otimes\overline{\mathfrak{g}}
\oplus\mathbb{C}K\oplus\mathbb{C}d
\end{equation*}
with Lie bracket
\[
[t^i\otimes x+\lambda K+\mu d,\;
 t^j\otimes y+\lambda'K+\mu'd]
\]
\begin{equation*}
= t^{i+j}\otimes[x,y]_0
+\mu j\,t^j\otimes y
-\mu' i\,t^i\otimes x
+i\delta_{i+j,0}(x|y)_0K,
\end{equation*}
for $i,j\in\mathbb{Z}$, $x,y\in\overline{\mathfrak{g}}$.

Let $\overline{\mathfrak{h}}$ be a Cartan subalgebra of $\overline{\mathfrak{g}}$
and $\overline{\Phi}$ its root system. For $\alpha\in\overline{\Phi}$ set
\begin{equation*}
\overline{\mathfrak{g}}_\alpha
=\{x\in\overline{\mathfrak{g}}\mid[h,x]_0=\alpha(h)x
\text{ for all }h\in\overline{\mathfrak{h}}\}.
\end{equation*}
Let $\overline{\Pi}=\{\overline{\alpha}_1,\dots,\overline{\alpha}_l\}$ and
$\overline{\Pi}^\vee=\{\overline{h}_1,\dots,\overline{h}_l\}$ be the sets of simple roots
and coroots. Then $\dim\overline{\mathfrak{g}}_\alpha=1$ for all
$\alpha\in\overline{\Phi}$. For each $i \in \{1,\dots,l\}$, let $\overline{e}_i$ be a basis vector of $\overline{\mathfrak{g}}_{\overline{\alpha}_i}$ and
$\overline{f}_i$ a basis vector of $\overline{\mathfrak{g}}_{-\overline{\alpha}_i}$. Then $\overline{\mathfrak{g}}$ is generated by
\[
\{\overline{h}_1,\dots,\overline{h}_l,\;
  \overline{e}_1,\dots,\overline{e}_l,\;
  \overline{f}_1,\dots,\overline{f}_l\}.
\]

Assume that $\overline{\mathfrak{g}}$ corresponds to a finite-type Cartan matrix
$\overline{A}$. Any permutation $\sigma\in S_l$ with
$\overline{a}_{ij}=\overline{a}_{\sigma(i)\sigma(j)}$ defines an automorphism by
\begin{equation*}
\overline{e}_i\mapsto\overline{e}_{\sigma(i)},\quad
\overline{f}_i\mapsto\overline{f}_{\sigma(i)},\quad
\overline{h}_i\mapsto\overline{h}_{\sigma(i)}.
\end{equation*}
Let $m=\mathrm{ord}(\sigma)$ and $\eta=e^{2\pi i/m}$. Define
$\tau\in\mathrm{Aut}(\hat{\overline{\mathfrak{g}}})$ by
\begin{equation*}
\tau(t^j\otimes x)=\eta^{-j}t^j\otimes\sigma(x),\quad
\tau(K)=K,\quad\tau(d)=d,
\end{equation*}
the \black{twisted automorphism}.

Let $\mathfrak{g}$ be an affine Kac--Moody algebra of type $X_N^{(r)}$ and
$\overline{\mathfrak{g}}$ of type $X_N$. If $r=1$, then
\begin{equation}\label{untwisted}
\mathfrak{g}\simeq\hat{\overline{\mathfrak{g}}}.
\end{equation}
If $r=2,3$, then $\sigma$ has order $r$ and
\begin{equation}\label{twisted}
\mathfrak{g}\simeq\hat{\overline{\mathfrak{g}}}^{\langle\tau\rangle}.
\end{equation}

The simple coroots $h_1,\dots,h_l$ of $\mathfrak{g}$ satisfy
\begin{equation}\label{hiplus0}
h_i\in1\otimes\overline{\mathfrak{h}},\qquad i=1,\dots,l.
\end{equation}
For details see \cite[Theorems~18.5, 18.9, 18.14]{Carter}.

\subsection{Dominant integral weights}

Define the \black{set of integral weights} by
\[
P=\sum_{i\in I}\mathbb{Z}\Lambda_i+\mathbb{C}\delta.
\]
For any $S\subset\mathfrak{h}^*$, let $\overline{S}$ be the subset of
$S+\mathbb{C}\delta$ consisting of elements $\lambda$ with $\lambda(d)=0$.
Then
\begin{equation*}\label{Pplusbar}
\overline{P}=\sum_{i\in I}\mathbb{Z}\Lambda_i.
\end{equation*}
For $\lambda\in P$, the integer $\lambda(K)$ is called the \black{level} of
$\lambda$. For $m\in\mathbb{Z}$, let $P^m$ denote the set of integral weights of
level $m$, that is,
\[
P^m
=\left\{\sum_{i\in I} m_i\Lambda_i \,\middle|\,
\sum_{i\in I} m_i c_i=m,\; m_i\in\mathbb{Z}\right\}
+\mathbb{C}\delta.
\]

Define the \black{set of dominant integral weights} by
\[
P_+=\sum_{i\in I}\mathbb{Z}_{\ge0}\Lambda_i+\mathbb{C}\delta.
\]
Let $\overline{P_+}$ and $\overline{P^m_+}$ be the intersections of $P_+$ with
$\overline{P}$ and $\overline{P^m}$, respectively.
\subsection{Highest weight representations of affine Kac--Moody algebras}

Let $\mathfrak{g}$ be an affine Kac--Moody algebra.  
The category $\mathcal{O}$ consists of $\mathfrak{g}$-modules $V$ such that
\begin{enumerate}
\item $V=\bigoplus_{\lambda\in\mathfrak{h}^*} V_\lambda$, where
\[
V_\lambda=\{v\in V \mid hv=\lambda(h)v \text{ for all } h\in\mathfrak{h}\};
\]
\item $\dim V_\lambda<\infty$ for all $\lambda\in\mathfrak{h}^*$;
\item there exist $\lambda_1,\dots,\lambda_s\in\mathfrak{h}^*$ such that
$V_\lambda\neq0$ implies $\lambda\le\lambda_i$ for some $i$.
\end{enumerate}
Morphisms in $\mathcal{O}$ are $\mathfrak{g}$-module homomorphisms.

For $\Lambda\in\mathfrak{h}^*$, the \black{Verma module} of highest weight
$\Lambda$ is
\[
M(\Lambda)=\mathcal{U}(\mathfrak{g})/K_\Lambda,
\]
where $\mathcal{U}(\mathfrak{g})$ is the \black{universal enveloping algebra} and
\[
K_\Lambda=\mathcal{U}(\mathfrak{g})\,\mathfrak{n}_+
+\sum_{h\in\mathfrak{h}}\mathcal{U}(\mathfrak{g})(h-\Lambda(h)).
\]
The module $M(\Lambda)$ has a unique maximal submodule $J(\Lambda)$, and we set
\[
L(\Lambda)=M(\Lambda)/J(\Lambda).
\]

\begin{proposition}
For all $\Lambda\in\mathfrak{h}^*$, one has $M(\Lambda)\in\mathcal{O}$ and hence
$L(\Lambda)\in\mathcal{O}$. Moreover, the $L(\Lambda)$ are exactly the irreducible
objects of $\mathcal{O}$.
\end{proposition}

A $\mathfrak{g}$-module $V$ is called \black{integrable} if
\[
V=\bigoplus_{\lambda\in\mathfrak{h}^*} V_\lambda
\]
and the operators $e_i$ and $f_i$ are locally nilpotent for all $i=1,\dots,l$.

\begin{proposition}
Let $\mathfrak{g}$ be an affine Kac--Moody algebra.  
Then $L(\Lambda)$ is integrable if and only if $\Lambda\in P_+$.
\end{proposition}
\subsection{Contravariant Hermitian forms on Verma modules}\label{w0}

Let $e_1,\dots,e_n,f_1,\dots,f_n,h_1,\dots,h_n,d$ be the \black{Chevalley generators} of $\mathfrak{g}$.  
Define the \black{$\mathbb{C}$-antilinear anti-involution} $\omega_0$ of $\mathfrak{g}$ by
\[
\omega_0(e_i)=f_i,\qquad
\omega_0(f_i)=e_i,\qquad
\omega_0(h_i)=h_i,\qquad
\omega_0(d)=d.
\]
Thus $\omega_0([x,y])=[\omega_0(y),\omega_0(x)]$ and $\omega_0^2=\mathrm{Id}$.  
The map $\omega_0$ extends uniquely to an antilinear anti-involution of
$\mathcal{U}(\mathfrak{g})$, again denoted by $\omega_0$.

By the triangular decomposition
$\mathfrak{g}=\mathfrak{n}_-\oplus\mathfrak{h}\oplus\mathfrak{n}_+$ and the
Poincaré--Birkhoff--Witt theorem,
\begin{equation}\label{decom-wo}
\mathcal{U}(\mathfrak{g})
=
\mathcal{U}(\mathfrak{h})
\oplus
\bigl(\mathfrak{n}_-\,\mathcal{U}(\mathfrak{g})
+
\mathcal{U}(\mathfrak{g})\,\mathfrak{n}_+\bigr).
\end{equation}
This decomposition is $\omega_0$-stable: $\omega_0$ fixes
$\mathcal{U}(\mathfrak{h})$ and interchanges
$\mathfrak{n}_-\,\mathcal{U}(\mathfrak{g})$ with
$\mathcal{U}(\mathfrak{g})\,\mathfrak{n}_+$. Since $\mathfrak{h}$ is abelian,
$\mathcal{U}(\mathfrak{h})=\mathcal{S}(\mathfrak{h})$.

Let $\pi\colon \mathcal{U}(\mathfrak{g})\to\mathcal{S}(\mathfrak{h})$ be the
projection onto the first summand in \eqref{decom-wo}. It is an
$\mathcal{S}(\mathfrak{h})$-bimodule map, i.e.,
\[
\pi(xgy)=x\,\pi(g)\,y
\qquad
(x,y\in\mathcal{S}(\mathfrak{h}),\ g\in\mathcal{U}(\mathfrak{g})).
\]
Define
\[
\langle x,y\rangle=\pi\bigl(\omega_0(x)y\bigr),
\qquad
x,y\in\mathcal{U}(\mathfrak{g}).
\]

\begin{proposition}\label{lrw0}
The form $\langle\cdot,\cdot\rangle$ satisfies:
\begin{enumerate}
\item[1.] antilinearity in the first variable and linearity in the second;
\item[2.] $\langle gx,y\rangle=\langle x,\omega_0(g)y\rangle$ for all $x,y,g\in\mathcal{U}(\mathfrak{g})$;
\item[3.] $\langle x,y\rangle=0$ if $x\in\mathcal{U}(\mathfrak{g})\mathfrak{n}_+$ or $y\in\mathcal{U}(\mathfrak{g})\mathfrak{n}_+$.
\end{enumerate}
\end{proposition}

Let $M$ be a $\mathfrak{g}$-module. A \black{Hermitian contravariant form}
with respect to $\omega_0$ is a map
\[
H\colon M\times M\to\mathbb{C}
\]
that is antilinear in the first variable, linear in the second, and satisfies
\[
H(gx,y)=H(x,\omega_0(g)y)
\quad (x,y\in M,\ g\in\mathfrak{g}).
\]

For $\Lambda\in\mathfrak{h}^*$, let
$ev_\Lambda\colon\mathcal{S}(\mathfrak{h})\to\mathbb{C}$ be defined by
\[
ev_\Lambda(x_1\cdots x_r)=\Lambda(x_1)\cdots\Lambda(x_r),
\qquad x_i\in\mathfrak{h}.
\]
Its kernel is generated by $h-\Lambda(h)$, $h\in\mathfrak{h}$. Set
\[
\langle\cdot,\cdot\rangle_\Lambda
=
ev_\Lambda\circ\langle\cdot,\cdot\rangle.
\]

Let $\mathfrak{h}_{\mathbb{R}}^*$ be the dual of
$\mathbb{R}d\oplus\bigoplus_{i=1}^n\mathbb{R}h_i$.
Using Proposition~\ref{lrw0}, we obtain the following result.

\begin{proposition}\label{hR*}
We have $\langle x,y\rangle_\Lambda=0$ for all $x$ or $y$ in
\[
\mathcal{U}(\mathfrak{g})\mathfrak{n}_+
+
\sum_{h\in\mathfrak{h}}\mathcal{U}(\mathfrak{g})(h-\Lambda(h))
\]
if and only if $\Lambda\in\mathfrak{h}_{\mathbb{R}}^*$.
\end{proposition}

Assume $\Lambda\in\mathfrak{h}_{\mathbb{R}}^*$. Then $\langle\cdot,\cdot\rangle_\Lambda$
descends to a Hermitian form on the Verma module $M(\Lambda)$, denoted again by
$\langle\cdot,\cdot\rangle_\Lambda$.

\begin{proposition}\label{lrLambdaM}
Let $\Lambda\in\mathfrak{h}_{\mathbb{R}}^*$. Then:
\begin{enumerate}
\item[1.] $\langle\cdot,\cdot\rangle_\Lambda$ is a Hermitian contravariant form on $M(\Lambda)$.
\item[2.] $\langle x,y\rangle_\Lambda=0$ for $x\in M(\Lambda)_{\Lambda-\alpha}$ and
$y\in M(\Lambda)_{\Lambda-\beta}$ with $\alpha\neq\beta$ in $Q_+$.
\item[3.] The radical of $\langle\cdot,\cdot\rangle_\Lambda$ equals the maximal submodule
$J(\Lambda)$.
\end{enumerate}
\end{proposition}

A $\mathfrak{g}$-module is \black{unitarizable} if it admits a positive definite
Hermitian form contravariant with respect to $\omega_0$.
Let $\Lambda\in\mathfrak{h}_{\mathbb{R}}^*$. By Proposition~\ref{lrLambdaM}(3),
$\langle\cdot,\cdot\rangle_\Lambda$ induces such a form on $L(\Lambda)$.
The following criterion is due to \cite[Theorem~11.7]{Kac}.

\begin{theorem}\label{positiveform}
Let $\Lambda\in\mathfrak{h}_{\mathbb{R}}^*$. The module $L(\Lambda)$ is unitarizable
if and only if $\Lambda\in P_+$.
\end{theorem}


\section{Branching on Cartan subalgebras}\label{Sect3}

We recall basic facts on branching with respect to Cartan subalgebras of affine
Kac--Moody algebras.

Let $\Lambda\in P_+$ and let $\mathfrak{h}$ be a Cartan subalgebra of $\mathfrak{g}$.
As an $\mathfrak{h}$-module, $L(\Lambda)$ decomposes into weight spaces
\begin{equation*}\label{Cartandecomposition}
L(\Lambda)=\bigoplus_{\lambda\in\mathfrak{h}^*} L(\Lambda)_\lambda,
\end{equation*}
where
\[
L(\Lambda)_\lambda
=
\{\,v\in L(\Lambda)\mid hv=\lambda(h)v \text{ for all } h\in\mathfrak{h}\,\}.
\]
Set $\mult_{\Lambda,\mathfrak{h}}(\lambda)=\dim L(\Lambda)_\lambda$.
The \black{character} of $L(\Lambda)$ is
\[
\ch_\Lambda=\sum_{\lambda\in\mathfrak{h}^*}
\mult_{\Lambda,\mathfrak{h}}(\lambda)\,e^\lambda .
\]
The \black{set of weights} is
\begin{equation*}\label{Cartanweights}
P(\Lambda)=\{\lambda\in\mathfrak{h}^*\mid \mult_{\Lambda,\mathfrak{h}}(\lambda)\neq0\}.
\end{equation*}

A weight $\lambda\in P(\Lambda)$ is \black{maximal} if
$\lambda+n\delta\notin P(\Lambda)$ for all $n>0$.
Let $\max(\Lambda)$ denote the set of maximal weights. Then
\begin{equation}\label{maxLamda}
\max(\Lambda)=W\bigl(\max(\Lambda)\cap P_+\bigr),
\end{equation}
and
\begin{equation}\label{PLamda}
P(\Lambda)
=
W\bigl((\Lambda-Q_+)\cap P_+\bigr)
=
\max(\Lambda)-\mathbb{Z}_{\ge0}\delta .
\end{equation}

For $\lambda\in P(\Lambda)+\mathbb{C}\delta$, let $b_{\Lambda,\lambda}\in\mathbb{C}$
be defined by $\lambda+b_{\Lambda,\lambda}\delta\in\max(\Lambda)$.
If $\lambda\in\Lambda+Q$, then $b_{\Lambda,\lambda}\in\mathbb{Z}$.
Moreover, $\lambda+b\delta\in P(\Lambda+b\delta)$ if and only if
$\lambda\in P(\Lambda)$, and
\begin{equation}\label{bmodule}
b_{\Lambda+b_2\delta,\lambda+b_1\delta}
=
b_{\Lambda,\lambda}+b_2-b_1
\end{equation}
for all $b_1,b_2\in\mathbb{C}$.

Let $\Gamma(\mathfrak{g},\mathfrak{h})$ be the set of pairs
$(\Lambda,\lambda)\in P_+\times\mathfrak{h}^*$ with $\lambda\in P(\Lambda)$.
\subsection{About the character $\ch_\Lambda$}

We recall basic properties of the character $\ch_\Lambda$ of an affine
Kac--Moody algebra.

Set $q=e^{-\delta}$. For $\lambda\in\mathfrak{h}^*$, define the
\black{string function}
$c^\Lambda_\lambda\in\mathbb{C}((q))$ by
\[
c^\Lambda_\lambda
=
\sum_{n\in\mathbb{Z}}
\mult_{\Lambda,\mathfrak{h}}(\lambda-n\delta)\,q^n .
\]
For all $w\in W$,
\begin{equation}\label{stringinv}
c^\Lambda_\lambda=c^\Lambda_{w\lambda},
\end{equation}
and the character decomposes as
\begin{equation}\label{charstring}
\ch_\Lambda
=
\sum_{\lambda\in\max(\Lambda)} c^\Lambda_\lambda\,e^\lambda .
\end{equation}

For $w\in W$, let $l(w)$ be the length of $w$ with respect to
$s_0,\dots,s_l$, and set $\epsilon(w)=(-1)^{l(w)}$.
The \black{Weyl--Kac character formula} gives
\begin{equation}\label{charkac}
\ch_\Lambda
=
\frac{\sum_{w\in W}\epsilon(w)\,e^{w(\Lambda+\rho)}}
     {\sum_{w\in W}\epsilon(w)\,e^{w(\rho)}} .
\end{equation}
\subsection{Semigroup structure}

We study the set $\Gamma(\mathfrak{g},\mathfrak{h}) \subset \mathfrak{h}^*\times\mathfrak{h}^*$.

\begin{theorem}\label{semigroupAPA}
The set $\Gamma(\mathfrak{g},\mathfrak{h})$ is a subsemigroup of
$\mathfrak{h}^*\times\mathfrak{h}^*$.
\end{theorem}

\begin{proof}
Let $(\Lambda,\lambda),(\overline{\Lambda},\overline{\lambda})\in
\Gamma(\mathfrak{g},\mathfrak{h})$.
Then $\lambda+\overline{\lambda}$ is a weight of
$L(\Lambda)\otimes L(\overline{\Lambda})$.
Hence it is a weight of $L(\Lambda'')$ for some
\[
\Lambda'' \in \bigl((\Lambda+\overline{\Lambda})-Q_+\bigr)\cap P_+ .
\]
By~\eqref{PLamda}, $P(\Lambda'')\subset P(\Lambda+\overline{\Lambda})$.
Thus $\lambda+\overline{\lambda}\in P(\Lambda+\overline{\Lambda})$, and
$(\Lambda+\overline{\Lambda},\lambda+\overline{\lambda})\in
\Gamma(\mathfrak{g},\mathfrak{h})$.
\end{proof}

\begin{remark}
The same argument applies to any symmetrizable Kac--Moody algebra
$\mathfrak{g}$.
\end{remark}


\section{Branching on winding subalgebras: the general case}\label{Sect4}

In this section, we study the branching problem for winding subalgebras.
\subsection{Winding subalgebras of affine Kac--Moody algebras}\label{winding}

We recall the notion of winding subalgebras following~\cite{KacWaki}, which
play a central role in tensor product decompositions and in the solution
of Frenkel’s conjecture.

Let $\mathfrak{g}$ be an affine Kac--Moody algebra of type $X_N^{(r)}$,
defined by~\eqref{untwisted} and~\eqref{twisted}.
Fix $u\in\mathbb{Z}_{>1}$ with $u\equiv 1 \pmod r$.
The assignment
\[
t^j\otimes x \mapsto t^{uj}\otimes x,\qquad
K \mapsto uK,\qquad
d \mapsto \frac{d}{u},
\]
for $j\in\mathbb{Z}$ and $x\in\overline{\mathfrak{g}}$, defines an injective Lie algebra homomorphism \(\psi_u\colon \mathfrak{g}\rightarrow \mathfrak{g}\).

For $r>1$, its image is stable under $\langle\tau\rangle$.

Define a subalgebra $\mathfrak{g}[u]\subset\mathfrak{g}$ by
\begin{enumerate}
\item if $r=1$, $\mathfrak{g}[u]=\psi_u(\hat{\overline{\mathfrak{g}}})$;
\item if $r>1$, $\mathfrak{g}[u]=\psi_u(\hat{\overline{\mathfrak{g}}})^{\langle\tau\rangle}$.
\end{enumerate}
Then $\mathfrak{g}[u]$ is isomorphic to $\mathfrak{g}$.
We call $\mathfrak{g}[u]$ the \black{winding subalgebra} of $\mathfrak{g}$
associated with $u$.
\subsection{Formulas for characteristic elements}

Set $\dot K=\psi_u(K)=uK$.  
Let $\tilde\psi_u\colon\mathfrak{h}\to\mathfrak{h}$ be the restriction of $\psi_u$.
For $i\in I$, set $\dot h_i=\tilde\psi_u(h_i)$. Then, by~\eqref{hiplus0},
\begin{equation*}\label{hdot}
\dot h_i=h_i \ \text{for } i>0, 
\qquad 
\dot h_0=\frac{u-1}{c_0}K+h_0 .
\end{equation*}

Let ${}^t\tilde\psi_u\colon\mathfrak{h}^*\to\mathfrak{h}^*$ be the dual map,
\begin{equation*}\label{dualmap}
{}^t\tilde\psi_u(\lambda)(h)=\lambda\bigl(\tilde\psi_u(h)\bigr),
\qquad \lambda\in\mathfrak{h}^*,\ h\in\mathfrak{h}.
\end{equation*}
Set $\dot\alpha_i={}^t\tilde\psi_u(\alpha_i)$. Then
\[
\dot\alpha_i=\alpha_i \ \text{for } i>0,
\qquad
\dot\alpha_0=\frac{u-1}{a_0}\delta+\alpha_0 .
\]

Define $\dot\Lambda_i={}^t\tilde\psi_u(\Lambda_i)$ and
$\dot\rho={}^t\tilde\psi_u(\rho)$. We obtain
\[
\dot\Lambda_i=\Lambda_i+\Bigl(\frac1u-1\Bigr)\frac{c_i}{c_0}\Lambda_0,
\qquad
\dot\rho=\rho+\Bigl(\frac1u-1\Bigr)\frac{\mathrm{h}^\vee}{c_0}\Lambda_0 .
\]

The map ${}^t\tilde\psi_u$ defines reflections
\[
\dot s_i(\lambda)=\lambda-\lambda(\dot h_i)\dot\alpha_i .
\]
The Weyl group $\dot W$ of $\mathfrak{g}[u]$, generated by $\dot s_i$ $(i\in I)$,
satisfies
\[
\dot W \cong t_{uM}\rtimes\overline W ,
\]
and hence $\dot W\subset W$ by~\eqref{Weylgroup}.

Let
\[
\dot P_+=\sum_{i\in I}\mathbb{Z}_{\ge0}\dot\Lambda_i+\mathbb{C}\delta
\]
be the set of dominant integral weights of $\mathfrak{g}[u]$.
For $m\in\mathbb{Z}_{\ge0}$, set
\[
\dot P_+^m
=\Bigl\{\sum_{i\in I} m_i\dot\Lambda_i \ \Bigm|\ 
\sum_{i\in I} m_i c_i=m,\ m_i\in\mathbb{Z}_{\ge0}\Bigr\}
+\mathbb{C}\delta .
\]
For $\lambda\in\dot P_+$, denote by $\dot L(\lambda)$ the irreducible integrable
$\mathfrak{g}[u]$-module of highest weight $\lambda$.
The winding subalgebra admits the triangular decomposition
\[
\mathfrak{g}[u]=\dot{\mathfrak n}_-\oplus\mathfrak h\oplus\dot{\mathfrak n}_+ .
\]
\subsection{The set of weights $P_u(\Lambda)$}\label{weightsPAu}

Let $\Lambda\in P_+^m$ with $m\in\mathbb{Z}_{\geq 0}$. The $\mathfrak{g}$-module $L(\Lambda)$ is a $\mathfrak{g}[u]$-module of level $um$. It decomposes into irreducible integrable $\mathfrak{g}[u]$-modules as
\begin{equation*}\label{windingdecomposition}
L(\Lambda)=\bigoplus_{\lambda\in \dot{P}_+^{um}}
\dot{L}(\lambda)^{\mult_{\Lambda,\mathfrak{g}[u]}(\lambda)}.
\end{equation*}

Define
\[
P_u(\Lambda)=\{\lambda\in\dot{P}_+ \mid \mult_{\Lambda,\mathfrak{g}[u]}(\lambda)\neq 0\}.
\]
Then $P_u(\Lambda)\subset P(\Lambda)$ \cite{KacWaki}. Hence
\[
P_u(\Lambda)\subset \max(\Lambda)-\mathbb{Z}_{\geq 0}\delta.
\]

A weight $\lambda\in P_u(\Lambda)$ is called a \black{$\mathfrak{g}[u]$-maximal weight} of $\Lambda$ if there exists no $n\in\mathbb{Z}_{>0}$ such that $\lambda+n\delta\in P_u(\Lambda)$. Denote by $\max_u(\Lambda)$ the set of all $\mathfrak{g}[u]$-maximal weights.

For each $\lambda\in P_u(\Lambda)+\mathbb{C}\delta$, there exists a unique $b_{\Lambda,\lambda,u}\in\mathbb{C}$ such that
\[
\lambda+b_{\Lambda,\lambda,u}\delta\in \max_u(\Lambda).
\]
By definition,
\begin{equation}\label{bulessthanb}
b_{\Lambda,\lambda}-b_{\Lambda,\lambda,u}\in\mathbb{Z}_{\geq 0}.
\end{equation}

Finally, let $\Gamma(\mathfrak{g},\mathfrak{g}[u])$ be the set of all pairs
$(\Lambda,\lambda)\in P_+\times \dot{P}_+$ with $\lambda\in P_u(\Lambda)$.
\subsection{Character method}

We recall basic facts on the representation theory of the Virasoro algebra needed for the branching problem. The approach follows \cite{KumarBrown, KacWaki}. Our goal is to express the character of an irreducible highest weight $\mathfrak{g}$-module $L(\Lambda)$, $\Lambda\in P_+$, in terms of the characters of irreducible highest weight $\mathfrak{g}[u]$-modules $\dot{L}(\lambda)$, $\lambda\in\dot{P}_+$. Virasoro characters enter naturally and encode the branching rule.
\subsubsection{The Virasoro algebra}

The \black{Virasoro algebra} $\Vir$ is the complex Lie algebra generated by $\{L_n \mid n\in\mathbb{Z}\}$ and $Z$, with relations
\begin{equation*}
[L_m,L_n]=(m-n)L_{m+n}+\frac{m^3-m}{12}\delta_{m+n,0}Z,
\qquad [\Vir,Z]=0.
\end{equation*}
Set $\Vir_0=\mathbb{C}L_0\oplus\mathbb{C}Z$.

Let $V$ be a \black{$\Vir$-module}. For $\lambda\in \Vir_0^*$, define the \black{weight space}
\begin{equation*}
V_\lambda=\{v\in V\mid Xv=\lambda(X)v \text{ for all } X\in \Vir_0\}.
\end{equation*}
Assume that $V$ admits a weight decomposition with finite-dimensional weight spaces.

Let $\omega_0^{\Vir}$ be the \black{$\mathbb{C}$-anti-linear anti-involution} of $\Vir$ given by
\begin{equation*}
\omega_0^{\Vir}(L_n)=L_{-n}, \qquad \omega_0^{\Vir}(Z)=Z.
\end{equation*}
Thus $\omega_0^{\Vir}([X,Y])=[\omega_0^{\Vir}(Y),\omega_0^{\Vir}(X)]$ and $(\omega_0^{\Vir})^2=\mathrm{Id}$.

A Hermitian form $\langle\cdot,\cdot\rangle$ on $V$ is \black{contravariant} with respect to $\omega_0^{\Vir}$ if
\begin{equation*}
\langle Xv,w\rangle=\langle v,\omega_0^{\Vir}(X)w\rangle
\end{equation*}
for all $v,w\in V$ and $X\in \Vir$. The module $V$ is \black{unitarizable} if such a form exists and is positive definite. Any unitarizable $\Vir$-module is completely reducible.

A $\Vir$-module $V$ is a \black{highest weight module} if there exists $0\neq v_0\in V$ such that $v_0$ is a $\Vir_0$-eigenvector, $L_n v_0=0$ for all $n>0$, and
\[
\mathcal{U}\!\left(\bigoplus_{n<0}\mathbb{C}L_n\right)v_0=V.
\]
Its highest weight $\lambda\in \Vir_0^*$ is defined by $Xv_0=\lambda(X)v_0$ for all $X\in \Vir_0$.

Let $V$ be a unitarizable highest weight $\Vir$-module with highest weight $\lambda$. Then
\[
\lambda\in \Vir_{0,\mathbb{R}}^*=\mathbb{R}(L_0)^*\oplus\mathbb{R}Z^*,
\]
and $\lambda(L_0),\lambda(Z)\ge 0$, since
\begin{equation*}
0\le \langle L_{-n}v_0,L_{-n}v_0\rangle
=\left(2n\lambda(L_0)+\frac{n^3-n}{12}\lambda(Z)\right)\langle v_0,v_0\rangle .
\end{equation*}

Let $\{(L_0)^*,Z^*\}$ be the dual basis of $\{L_0,Z\}$. We recall \cite[Lemma~4.1]{KumarBrown}.

\begin{lemma}\label{virakumar}
Let $V$ be a unitarizable irreducible highest weight $\Vir$-module with highest weight $\lambda$. Then:
\begin{enumerate}
\item[1.] If $\lambda(L_0)\neq 0$, then $V_{\lambda+n(L_0)^*}\neq 0$ for all $n\in\mathbb{Z}_{\ge 0}$.
\item[2.] If $\lambda(L_0)=0$ and $\lambda(Z)\neq 0$, then $V_{\lambda+n(L_0)^*}\neq 0$ for all $n\in\mathbb{Z}_{>1}$ and $V_{\lambda+(L_0)^*}=0$.
\item[3.] If $\lambda(L_0)=\lambda(Z)=0$, then $V$ is one-dimensional.
\end{enumerate}
\end{lemma}
\subsubsection{Unitarizability}

Let $\Lambda\in P_+$ and $\lambda\in \dot{P}_+ + \mathbb{C}\delta$. Define
\begin{equation*}
\mathcal{U}(\Lambda,\lambda)
=\left\{ v\in L(\Lambda)\ \middle|\ 
(\mathfrak{n}_+\cap\mathfrak{g}'[u])v=0,\;
hv=\lambda(h)v \ \forall\, h\in \mathfrak{h}\cap\mathfrak{g}'[u] \right\}.
\end{equation*}
By Propositions \ref{Lu-untwisted} and \ref{Lu-twisted}, the operators $L_n^{[u]}$ commute with $\mathfrak{g}'[u]$; hence $\mathcal{U}(\Lambda,\lambda)$ is $\Vir$-stable. Moreover, $L(\Lambda)$ decomposes as a $\mathfrak{g}'[u]\oplus \Vir$-module \cite{KacWaki, KacWaki2}
\begin{equation*}
L(\Lambda)
=\bigoplus_{\lambda\in \dot{P}_+ \bmod \mathbb{C}\delta}
\dot{L}(\lambda)\otimes \mathcal{U}(\Lambda,\lambda).
\end{equation*}

The positive definite contravariant Hermitian form $\langle\cdot,\cdot\rangle_\Lambda$ from Subsection~\ref{w0} restricts to $\mathcal{U}(\Lambda,\lambda)$. By Proposition~\ref{w0LZuu}, it is contravariant with respect to $\Vir$:
\begin{equation*}
\langle L_n^{[u]}v,w\rangle_\Lambda
=\langle v,L_{-n}^{[u]}w\rangle_\Lambda
=\langle v,\omega_0^{\Vir}(L_n)w\rangle_\Lambda,
\end{equation*}
\begin{equation*}
\langle Z^{[u]}v,w\rangle_\Lambda
=\langle v,Z^{[u]}w\rangle_\Lambda
=\langle v,\omega_0^{\Vir}(Z)w\rangle_\Lambda.
\end{equation*}
Therefore, $\mathcal{U}(\Lambda,\lambda)$ is a unitarizable $\Vir$-module.
\subsubsection{An identity of characters}

Let $\Lambda\in P_+^m$ with $m\in\mathbb{Z}_{\geq 0}$. Using
\eqref{maxLamda}, \eqref{stringinv}, \eqref{charstring}, and the inclusion
$\dot{W}\subset W$, we obtain
\begin{equation}\label{char.Ap}
\Bigl(\sum_{w\in\dot{W}}\epsilon(w)e^{w(\dot{\rho})}\Bigr) \ch_\Lambda
=
\sum_{\lambda\in\max(\Lambda)}
\Bigl(\sum_{w\in\dot{W}}\epsilon(w)e^{w(\lambda+\dot{\rho})}\Bigr)
c^\Lambda_\lambda .
\end{equation}

Assume that $\lambda+\dot{\rho}$ is $\dot{W}$-regular. Then there exist unique
$\sigma\in\dot{W}$ and $\lambda'\in\dot{P}_+$ such that
\[
\sigma(\lambda+\dot{\rho})=\lambda'+\dot{\rho}.
\]
Set $p(\lambda)=\epsilon(\sigma)$ and $\{\lambda\}=\lambda'$. If
$\lambda+\dot{\rho}$ is nonregular, set $p(\lambda)=0$. Then
\begin{equation*}
\sum_{w\in\dot{W}}\epsilon(w)e^{w(\lambda+\dot{\rho})}
=
p(\lambda)\sum_{w\in\dot{W}}\epsilon(w)e^{w(\{\lambda\}+\dot{\rho})}.
\end{equation*}
Combining \eqref{charkac} and \eqref{char.Ap} yields:

\begin{proposition}\label{propchar}
\begin{equation}\label{char.p}
\ch_\Lambda
=
\sum_{\lambda\in\max(\Lambda)}
p(\lambda)\,\dot{\ch}_{\{\lambda\}}\,c^\Lambda_\lambda .
\end{equation}
\end{proposition}
\subsection{Semigroup structure}

We state the first result on the set $\Gamma(\mathfrak{g},\mathfrak{g}[u])$.

\begin{theorem}\label{M1}
The set $\Gamma(\mathfrak{g},\mathfrak{g}[u])$ is a subsemigroup of $\mathfrak{h}^* \times \mathfrak{h}^*$.
\end{theorem}

\begin{proof}
Let $(\Lambda,\lambda),(\overline{\Lambda},\overline{\lambda})
\in\Gamma(\mathfrak{g},\mathfrak{g}[u])$. We show that
\[
(\Lambda+\overline{\Lambda},\lambda+\overline{\lambda})
\in\Gamma(\mathfrak{g},\mathfrak{g}[u]).
\]
By definition, $(\Lambda,\lambda)\in\Gamma(\mathfrak{g},\mathfrak{g}[u])$
if and only if $\dot{L}(\lambda)$ occurs in $L(\Lambda)$, equivalently,
there exists $0\neq v\in L(\Lambda)$ such that
\begin{equation}\label{cond.semi.u}
g v=0 \quad \forall g\in\dot{\mathfrak{n}}_+,
\qquad
h v=\lambda(h)v \quad \forall h\in\mathfrak{h}.
\end{equation}
Let $0\neq\overline{v}\in L(\overline{\Lambda})$ satisfy the analogous
conditions for $(\overline{\Lambda},\overline{\lambda})$. It suffices to
construct a nonzero vector $\tilde{v}\in L(\Lambda+\overline{\Lambda})$
satisfying \eqref{cond.semi.u} with weight
$\lambda+\overline{\lambda}$. We proceed in the following steps.

\noindent\textbf{Step 1. Construction of $\tilde{v}$.}
The module $L(\Lambda+\overline{\Lambda})$ occurs in
$L(\Lambda)\otimes L(\overline{\Lambda})$ with multiplicity one. Hence there
exists a unique $\mathfrak{g}$-stable complement $S$ such that
\[
L(\Lambda)\otimes L(\overline{\Lambda})
=
L(\Lambda+\overline{\Lambda}) \oplus S .
\]
Let
\[
\pi: L(\Lambda)\otimes L(\overline{\Lambda}) \to L(\Lambda+\overline{\Lambda})
\]
be the projection with kernel $S$, and define
\[
\tilde{v}=\pi(v\otimes\overline{v}).
\]
In the following steps we show that $\tilde{v}\neq 0$ and that it satisfies
\eqref{cond.semi.u}.

\noindent\textbf{Step 2. Nonvanishing of $\tilde{v}$.}
Write
\[
L(\Lambda)=\bigoplus_{\mu\in\mathfrak{h}^*}L(\Lambda)_\mu,
\qquad
L(\Lambda)^\vee=\bigoplus_{\mu\in\mathfrak{h}^*}(L(\Lambda)_\mu)^*.
\]
There exists $\psi\in L(\Lambda)^\vee\setminus\{0\}$ such that
\[
g\psi=0 \ \ (g\in\mathfrak{n}_-),\qquad
h\psi=-\Lambda(h)\psi \ \ (h\in\mathfrak{h}).
\]
Let $G$ be the minimal Kac--Moody group of $\mathfrak{g}$ \cite{Kumar2}. For $v\in L(\Lambda)$
define
\[
f_v(g)=\psi(g^{-1}v),\qquad g\in G.
\]
Since $L(\Lambda)$ is irreducible, $f_v\neq0$ (otherwise $Gv \subset \ker \psi$), and for $b\in B^-$,
\[
(1,b)\cdot f_v=\Lambda(b)^{-1}f_v.
\]

Similarly, choose $\overline{\psi}\in L(\overline{\Lambda})^\vee$ and define
$f_{\overline{v}}$. Then $f_{\overline{v}}\neq0$ and
\begin{equation*}\label{1bf}
(1,b)\cdot f_{\overline{v}}=\overline{\Lambda}(b)^{-1}f_{\overline{v}}
\qquad(b\in B^-).
\end{equation*}
Set $f=f_v f_{\overline{v}}$. Since $G$ is irreducible as an ind-variety, $f\neq0$, and
\begin{equation}\label{1bprodf}
(1,b)\cdot f=(\Lambda+\overline{\Lambda})(b)^{-1}f
\qquad(b\in B^-).
\end{equation}
Moreover,
\begin{equation}\label{prodf}
f(g)=(\psi\otimes\overline{\psi})\bigl(g^{-1}(v\otimes\overline{v})\bigr).
\end{equation}
We have
\[
L(\Lambda)^\vee\otimes L(\overline{\Lambda})^\vee
=(L(\Lambda)\otimes L(\overline{\Lambda}))^\vee
=L(\Lambda+\overline{\Lambda})^\vee\oplus S^\vee.
\]
By \eqref{1bprodf} and \eqref{prodf},
\[
\psi\otimes\overline{\psi}\in L(\Lambda+\overline{\Lambda})^\vee,
\qquad
\ker(\psi\otimes\overline{\psi})\supset S.
\]
Write $v \otimes \overline{v} = \pi(v \otimes \overline{v}) + s$ for some $s \in S$. Then
\begin{align*}
(\psi\otimes \overline{\psi})(g^{-1}(v\otimes \overline{v}))
&= (\psi\otimes \overline{\psi})(g^{-1}(\pi(v\otimes \overline{v})+s)) \\
&= (\psi\otimes \overline{\psi})(g^{-1}(\pi(v\otimes \overline{v}))).
\end{align*}
Hence $f(g) = (\psi\otimes \overline{\psi})(g^{-1}(\pi(v\otimes \overline{v})))$. Since $f\neq0$, we conclude $\tilde{v}=\pi(v\otimes\overline{v})\neq0$.

\noindent\textbf{Step 3. Verification of \eqref{cond.semi.u}.}
Let $g\in\dot{\mathfrak{n}}_+$ and $h\in\mathfrak{h}$. Since $\pi$ is a
$\mathfrak{g}$-module homomorphism, we have
\begin{align*}
g\,\tilde{v}
&= g\,\pi(v\otimes\overline{v})
 = \pi\bigl(g(v\otimes\overline{v})\bigr)
 = \pi\bigl(gv\otimes\overline{v}+v\otimes g\overline{v}\bigr)
 = 0, \\
h\,\tilde{v}
&= \pi\bigl(h(v\otimes\overline{v})\bigr)
 = \pi\bigl(hv\otimes\overline{v}+v\otimes h\overline{v}\bigr)
 = (\lambda+\overline{\lambda})(h)\,\tilde{v}.
\end{align*}
Thus $\tilde{v}$ satisfies \eqref{cond.semi.u}. This proves that
$\Gamma(\mathfrak{g},\mathfrak{g}[u])$ is a semigroup.
\end{proof}

\begin{remark}
The same argument applies to any subalgebra $\tilde{\mathfrak{g}}\subset\mathfrak{g}$
admitting a triangular decomposition
\[
\tilde{\mathfrak{g}}
=
(\tilde{\mathfrak{g}}\cap\mathfrak{n}_-)
\oplus
(\tilde{\mathfrak{g}}\cap\mathfrak{h})
\oplus
(\tilde{\mathfrak{g}}\cap\mathfrak{n}_+).
\]
\end{remark}
\subsection{Description of the set $\Gamma(\mathfrak{g},\mathfrak{g}[u])$}

Let $(\Lambda,\lambda)\in P_+\times\dot{P}_+$. Describing
$\Gamma(\mathfrak{g},\mathfrak{g}[u])$ is equivalent to determining
\[
\mathcal{B}(\Lambda,\lambda)
=
\{\, b\in\mathbb{C}\mid \dot{L}(\lambda+b\delta)\subset L(\Lambda)\,\}.
\]

We may assume $\Lambda\in\overline{P_+}$ and $\lambda\in\overline{\dot{P}_+}$.
Indeed,
\[
\dot{L}(\lambda+b_1\delta)\subset L(\Lambda+b_2\delta)
\iff
\dot{L}(\lambda+(b_1-b_2)\delta)\subset L(\Lambda),
\]
hence
\begin{equation}\label{BAa=BAa}
\mathcal{B}(\Lambda+b_2\delta,\lambda+b_1\delta)
=
\mathcal{B}(\Lambda,\lambda)+b_2-b_1 .
\end{equation}
Since $\mathcal{B}(\Lambda,\lambda)\neq\varnothing$ if and only if
$\lambda\in P_u(\Lambda)+\mathbb{C}\delta$, we further assume
$\lambda\in\overline{P_u(\Lambda)}$.

Let $m\in\mathbb{Z}_{>0}$ be the level of $\Lambda$. Set
\[
c_m=\frac{r m\,\dim\overline{\mathfrak{g}}}{m+\mathrm{h}^\vee},
\qquad
c_m^{[u]}=u c_m-c_{um},
\]
\[
m_\Lambda=\frac{|\Lambda+\rho|^2}{2(m+\mathrm{h}^\vee)}
-\frac{|\rho|^2}{2\mathrm{h}^\vee},
\qquad
\dot{m}_\lambda=\frac{|\lambda+\dot{\rho}|^2}{2(um+\mathrm{h}^\vee)}
-\frac{|\dot{\rho}|^2}{2\mathrm{h}^\vee},
\]
and
\begin{equation}\label{hAa}
h_{\Lambda,\lambda}^{[u]}
=
u^{-1}m_\Lambda-\dot{m}_\lambda+\frac{c_m^{[u]}}{24}-b_{\Lambda,\lambda,u}.
\end{equation}

The next theorem gives a description of $\mathcal{B}(\Lambda,\lambda)$.

\begin{theorem}\label{BAa}
Let $\mathfrak{g}$ be an affine Kac--Moody algebra of type $X_N^{(r)}$. Fix $u \in \mathbb{Z}_{>1}$ such that $u \equiv 1 \pmod r$. Let $\Lambda \in P_+$ and $\lambda \in P_u(\Lambda)+\mathbb{C}\delta$. Then:
\begin{enumerate}
    \item[1.] $b_{\Lambda,\lambda,u}-(\mathbb{Z}_{\geq 0}\setminus \{1\}) \subset \mathcal{B}(\Lambda,\lambda) \subset b_{\Lambda,\lambda,u}-\mathbb{Z}_{\geq 0}$.
    \item[2.] If, in addition, $h_{\Lambda,\lambda}^{[u]} \neq 0$, then $\mathcal{B}(\Lambda,\lambda) = b_{\Lambda,\lambda,u}-\mathbb{Z}_{\geq 0}$.
\end{enumerate}
\end{theorem}

\begin{proof}
By definition and since $P_u(\Lambda)\subset P(\Lambda)$, we have
\[
\mathcal{B}(\Lambda,\lambda)\subset b_{\Lambda,\lambda,u}-\mathbb{Z}_{\ge 0}.
\]
Using \eqref{BAa=BAa}, we may assume $\Lambda\in\overline{P_+^m}$ and
$\lambda\in\overline{P_u(\Lambda)}$. The $\mathfrak{g}'[u]\oplus \Vir$–module decomposition of $L(\Lambda)$ is
\[
L(\Lambda)
=\bigoplus_{\lambda\in\overline{P_u(\Lambda)}}
\dot{L}(\lambda)\otimes\mathcal{U}(\Lambda,\lambda),
\]
where $Z^{[u]}$ acts on $\mathcal{U}(\Lambda,\lambda)$ by scalar multiplication $c_m^{[u]}\neq 0$ and
the lowest $L_0^{[u]}$–eigenvalue equals $h_{\Lambda,\lambda}^{[u]}$
\cite{KacWaki}. By \cite{KacWaki, KacWaki2},
\begin{equation*}\label{charvir}
\ch_\Lambda
=\sum_{\lambda\in\overline{P_u(\Lambda)}}
\dot{\ch}_\lambda\,
q^{\dot m_\lambda-u^{-1}m_\Lambda}
\operatorname{tr}_{\mathcal{U}(\Lambda,\lambda)}
\!\left(q^{L_0^{[u]}-Z^{[u]}/24}\right).
\end{equation*}
Since the $\Vir$–-action on $\mathcal{U}(\Lambda,\lambda)$ is unitarizable, this
rewrites as
\begin{align*}
    \ch_\Lambda
    &= \sum_{\lambda \in \overline{P_u(\Lambda)}}
    \dot{\ch}_{\lambda}
    q^{\dot{m}_\lambda-u^{-1}m_\Lambda-c_m^{[u]}/24+h_{\Lambda,\lambda}^{[u]}}
    \bigl(\dim \mathcal{U}(\Lambda,\lambda)_{\lambda'}
    + \dim \mathcal{U}(\Lambda,\lambda)_{\lambda'+(L_0)^*} q + \cdots \bigr) \\
    &= \sum_{\lambda \in \overline{P_u(\Lambda)}}
    \dot{\ch}_{\lambda}
    q^{-b_{\Lambda,\lambda,u}}
    \bigl(\dim \mathcal{U}(\Lambda,\lambda)_{\lambda'}
    + \dim \mathcal{U}(\Lambda,\lambda)_{\lambda'+(L_0)^*} q + \cdots \bigr) \\
    &= \sum_{\lambda \in \overline{P_u(\Lambda)}}
    \dot{\ch}_{\lambda + b_{\Lambda,\lambda,u}\delta}
    \bigl(\dim \mathcal{U}(\Lambda,\lambda)_{\lambda'}
    + \dim \mathcal{U}(\Lambda,\lambda)_{\lambda'+(L_0)^*} q + \cdots \bigr).
\end{align*}

Let $V\subset\mathcal{U}(\Lambda,\lambda)$ be the highest weight $\Vir$–-module
with highest weight $\lambda'$ satisfying
\[
\lambda'(Z^{[u]})=c_m^{[u]},\qquad
\lambda'(L_0^{[u]})=h_{\Lambda,\lambda}^{[u]},
\]
which is contained in $\mathcal{U}(\Lambda,\lambda)$. By Lemma~\ref{virakumar}:
\begin{enumerate}
\item If $h_{\Lambda,\lambda}^{[u]}\neq 0$, then
\[
    0 \neq \dim V_{\lambda'+n(L_0^{[u]})^*}
    \leq \dim \mathcal{U}(\Lambda,\lambda)_{\lambda'+n(L_0^{[u]})^*}
    \quad \text{for all } n \in \mathbb{Z}_{\geq 0}.
    \]
Hence $\dot{L}(\lambda+(b_{\Lambda,\lambda,u}-n)\delta)\subset L(\Lambda)
\quad\text{for all }n\in\mathbb{Z}_{\ge 0}.$

\item If $h_{\Lambda,\lambda}^{[u]}=0$, the same holds for all
$n\in\mathbb{Z}_{\ge 0}\setminus\{1\}$.
\end{enumerate}
This completes the proof.
\end{proof}
\section{The cases $A^{(1)}_1$ and $A^{(2)}_2$}\label{A1A2}
In this section, we treat the cases $A^{(1)}_1$ and $A^{(2)}_2$.  
We compute $b_{\Lambda,\lambda}$ explicitly and determine when
$b_{\Lambda,\lambda,u}=b_{\Lambda,\lambda}$.  
This yields relations between $\Gamma(\mathfrak{g},\mathfrak{g}[u])$
and its saturation.

\subsubsection{Computation of $b_{\Lambda,\lambda}$ for $A^{(1)}_1$}\label{bAa.A11}

We compute $\max(\Lambda)$ and hence $b_{\Lambda,\lambda}$ for $A^{(1)}_1$.
By \eqref{bmodule}, we may assume $\Lambda \in \overline{P_+}$.
The computation follows \cite{KumarBrown}.

Let $\mathfrak{g}$ be of type $A^{(1)}_1$ and fix $m\in\mathbb{Z}_{>0}$.
Let $\alpha=\alpha_1$ be the simple root. Then
\begin{equation*}\label{Pm+.1.1}
\overline{P^m_+}
=\left\{\, m\Lambda_0+\frac{j\alpha}{2}\ \middle|\ j\in[0,m]\cap\mathbb{Z}\right\}.
\end{equation*}

The following description of $\max(\Lambda)$ and $b_{\Lambda,\lambda}$
combines \cite[Lemmas~5.2 and~5.3]{KumarBrown}.

\begin{proposition}\label{maxLamdaandnk1}
Let $\Lambda = m\Lambda_0 + \frac{j\alpha}{2} \in \overline{P^m_+}$. For each $k \in \mathbb{Z}$, let $\phi(m,j,k)$ be the number uniquely determined by $k,m,j$ as follows:
\begin{enumerate}
\item[1.] Write $k = mq + r$ for some $q \in \mathbb{Z}$ and $r \in [0,m)$.
\item[2.] Set
\begin{equation}\label{nk1}
\phi(m,j,k)
= -q(k+r+j) +
\begin{cases}
-r, & \text{if } r \in [0,m-j], \\
m-j-2r, & \text{if } r \in [m-j,m).
\end{cases}
\end{equation}
\end{enumerate}
Then
\begin{equation*}\label{maxLamda.1.1}
\max(\Lambda)
= \{ \Lambda + k\alpha + \phi(m,j,k)\delta \mid k \in \mathbb{Z} \}.
\end{equation*}
Equivalently, for each $\lambda = m\Lambda_0 + \frac{j'\alpha}{2}$ with $j' \in j + 2\mathbb{Z}$, we have
\begin{equation}\label{ba11}
b_{\Lambda,\lambda}
= \phi\!\left(m,j,\frac{j'-j}{2}\right).
\end{equation}
\end{proposition}

To prove the proposition, we use the following lemma. \cite[Proposition~4.4]{KumarBrown} states the result for untwisted affine Kac--Moody algebras.
The statement and proof extend to all affine Kac--Moody algebras.

\begin{lemma}\label{maxlem}
For any affine Kac--Moody algebra and any $\Lambda \in P_+$ of positive level,
\begin{equation*}\label{max}
\max(\Lambda)\cap P_+
=
\left\lbrace
\Lambda - \sum_{i\in I} m_i\alpha_i
\,\middle\vert\,
m_i \in \mathbb{Z}_{\geq 0} \text{ for all } i,
\; m_i < a_i \text{ for some } i \in I
\right\rbrace
\cap P_+.
\end{equation*}
\end{lemma}

Using Lemma~\ref{maxlem}, we prove Proposition~\ref{maxLamdaandnk1}.

\begin{proof}
We have $\max(\Lambda)=W(\max(\Lambda)\cap P_+)$. By Lemma~\ref{maxlem},
\begin{equation*}
\max(\Lambda)\cap P_+
=
\left\lbrace
\Lambda-m_0(\delta-\alpha),\;
\Lambda-m_1\alpha
\,\middle\vert\,
m_i\in\mathbb{Z}_{\geq 0},\;
m_0 \leq \frac{m-j}{2},\;
m_1 \leq \frac{j}{2}
\right\rbrace.
\end{equation*}
Recall that $W=\{t_{n\alpha},\, t_{n\alpha}s_1 \mid n \in \mathbb{Z}\}$. We compute
\begin{align}
t_{n\alpha}(\Lambda-m_0(\delta-\alpha))
&= \Lambda + (m_0+mn)\alpha - \bigl((j+2m_0+mn)n+m_0\bigr)\delta, \label{tna1} \\
t_{n\alpha}s_1(\Lambda-m_0(\delta-\alpha))
&= \Lambda + (-j-m_0+mn)\alpha - \bigl((-j-2m_0+mn)n+m_0\bigr)\delta, \label{tna2} \\
t_{n\alpha}(\Lambda-m_1\alpha)
&= \Lambda + (-m_1+mn)\alpha - (j-2m_1+mn)n\delta, \label{tna3} \\
t_{n\alpha}s_1(\Lambda-m_1\alpha)
&= \Lambda + (-j+m_1+mn)\alpha - (-j+2m_1+mn)n\delta. \label{tna4}
\end{align}
Thus any element $\lambda \in \max(\Lambda)$ has the form $\Lambda + r\alpha + n'\delta$ for some $r,n' \in \mathbb{Z}$.

Fix such a $\lambda$. For any $q \in \mathbb{Z}$,
\begin{equation*}
t_{q\alpha}(\Lambda+r\alpha+n'\delta)
=
\Lambda + (mq+r)\alpha + \bigl(n'-(j+2r+mq)q\bigr)\delta
\end{equation*}
is still in $\max(\Lambda)$. Writing $k=mq+r$, we obtain
\begin{equation*}
\Lambda + k\alpha + n''\delta \in \max(\Lambda),
\qquad
n'' = n' - q(k+r+j).
\end{equation*}
Assume now that $0 \leq r < m$, so that $k=mq+r$ is the Euclidean division. By \eqref{tna1}--\eqref{tna4}, we have
\begin{equation*}
n'
=
\begin{cases}
-r, & \text{if } r \in [0,m-j], \\
m-j-2r, & \text{if } r \in [m-j,m).
\end{cases}
\end{equation*}
Hence $n''=\phi(m,j,k)$ as defined in \eqref{nk1}. This proves that
\begin{equation*}
\max(\Lambda)
=
\{ \Lambda + k\alpha + \phi(m,j,k)\delta \mid k \in \mathbb{Z} \}.
\end{equation*}
\end{proof}
\subsubsection{The conditions for $b_{\Lambda,\lambda,u} = b_{\Lambda,\lambda}$ in the case $A^{(1)}_1$}

As we have seen in the previous part, we know explicitly the number $b_{\Lambda,\lambda}$ in the case $A^{(1)}_1$. By Theorem~\ref{BAa}, the number $b_{\Lambda,\lambda,u}$ is the key to understanding the support $\Gamma(\mathfrak{g},\mathfrak{g}[u])$. The explicit formula for the number $b_{\Lambda,\lambda,u}$ is not known in general, even in the case $A^{(1)}_1$. However, in this part, we give some conditions under which we know that $b_{\Lambda,\lambda,u}=b_{\Lambda,\lambda}$.

\begin{theorem}\label{b=b.a11}
Let $\mathfrak{g}$ be the affine Kac--Moody algebra of type $A^{(1)}_1$. Let $\Lambda = m\Lambda_0 + \frac{j\alpha}{2} \in P_+$ and let $\lambda = m'\Lambda_0 + \frac{j'\alpha}{2} \in \dot{P}_+$.
\begin{enumerate}
    \item[1.] If there exists $b\in \mathbb{C}$ such that $(\Lambda,\lambda+b\delta)$ belongs to $\Gamma(\mathfrak{g},\mathfrak{g}[u])$, then $j'-j \in 2\mathbb{Z}$ and $m'=m$.
    \item[2.] If moreover
    \begin{enumerate}
        \item[a.] $j \leq j' \leq um-j$ and $u$ is even; or
        \item[b.] $j \leq j' \leq um - (m-j)$ and $u$ is odd,
    \end{enumerate}
    then $b_{\Lambda,\lambda,u} = b_{\Lambda,\lambda}$.
\end{enumerate}
\end{theorem}

For $\Lambda$ as in Theorem~\ref{b=b.a11}, let $\mathcal{A}_u(\Lambda)$ be the set of $\lambda$ to which the theorem applies, namely
\begin{equation}\label{AuA.11}
\mathcal{A}_u(\Lambda)
=\left\{\, m\Lambda_0+\frac{j'\alpha}{2}\,\middle|\, j'\in [j,um-j^*]\cap (j+2\mathbb{Z}) \right\},
\end{equation}
where $j^*=j$ if $u$ is even and $j^*=m-j$ if $u$ is odd.

\begin{remark}
This result extends \cite[Theorem~2.2]{KacWaki} for $A^{(1)}_1$ from level one to arbitrary positive level.
\end{remark}

We now state auxiliary propositions.

\begin{proposition}\label{prop1}
Let $\mathfrak{g}$ be the affine Kac--Moody algebra of type $A^{(1)}_1$. Fix $m \in \mathbb{Z}_{>0}$ and $u\in \mathbb{Z}_{>1}$. Let
\[
\Lambda = m\Lambda_0 + \frac{j\alpha}{2} \in \overline{P^m_+}.
\]
\begin{enumerate}
\item[1.]
We parametrize $\lambda \in \max(\Lambda)$ such that $\lambda+\dot{\rho}$ is regular with respect to $\dot{W}$ by
\[
\lambda_k = \Lambda + k\alpha + \phi(m,j,k)\delta.
\]
Then the only possible values of $k$ are
\begin{equation*}
k=\frac{j'-j}{2}-n(um+2)
\quad\text{and}\quad
k=-\frac{j'+j}{2}-1+n(um+2),
\end{equation*}
where $j'\in [0,um]\cap (j+2\mathbb{Z})$ and $n\in \mathbb{Z}$.
\item[2.]
Let
\begin{equation*}
N_k=-\phi(m,j,k)+un(j'+1-num-2n).
\end{equation*}
Then
\begin{enumerate}
\item[(i)] If $k=\frac{j'-j}{2}-n(um+2)$, then $p(\lambda_k)=1$ and
\[
\{\lambda_k\}=m\Lambda_0+\frac{j'\alpha}{2}-N_k\delta.
\]
\item[(ii)] If $k=-\frac{j'+j}{2}-1+n(um+2)$, then $p(\lambda_k)=-1$ and
\[
\{\lambda_k\}=m\Lambda_0+\frac{j'\alpha}{2}-N_k\delta.
\]
\item[(iii)]
Viewed as a function of $n$, the function $N_k$ attains its minimum at $n=0$ in the first case, and at $n=0$ or $n=1$ in the second case.
\end{enumerate}
\end{enumerate}
\end{proposition}
\begin{proof}
We need the following data:
\begin{equation*}
|\alpha|^2=2, \quad 
\dot{\rho}=\frac{2}{u}\Lambda_0+\frac{1}{2}\alpha, \quad 
\dot{W}=\{t_{un\alpha},\, t_{un\alpha}s_1 \mid n\in \mathbb{Z}\}.
\end{equation*}
Since $\lambda_k+\dot{\rho}$ is regular with respect to $\dot{W}$, there exists a unique $\sigma \in \dot{W}$ and
\[
\mu=m\Lambda_0+\frac{j'\alpha}{2}+b'\delta \in \dot{P}^{um}_+
\]
such that $\sigma(\lambda_k+\dot{\rho})=\mu + \dot{\rho}$.

\begin{enumerate}
\item[(a)] (Proof of (1) and (2.i).)
If $\sigma =t_{un\alpha}$ for some $n\in \mathbb{Z}$, then
$\sigma(\lambda_k+\dot{\rho})-\dot{\rho}$ equals
\begin{equation*}\label{j'A11}
m\Lambda_0+\left(num+2n+k+\frac{j}{2}\right)\alpha
+\bigl(\phi(m,j,k)-un(2k+j+1+num+2n)\bigr)\delta.
\end{equation*}
Hence
\begin{equation*}
j'\in [0,um]\cap (j+2\mathbb{Z})
\quad\text{and}\quad
k=\frac{j'-j}{2}-n(um+2).
\end{equation*}
In this case,
\begin{equation*}
p(\lambda_k)=1
\quad\text{and}\quad
\{\lambda_k\}
=m\Lambda_0+\frac{j'\alpha}{2}
+\bigl(\phi(m,j,k)-un(j'+1-num-2n)\bigr)\delta.
\end{equation*}

\item[(b)] (Proof of (1) and (2.ii).)
If $\sigma =t_{un\alpha}s_1$ for some $n\in \mathbb{Z}$, then
$\sigma(\lambda_k+\dot{\rho})-\dot{\rho}$ equals
\begin{equation*}
m\Lambda_0+\left(num+2n-k-\frac{j}{2}-1\right)\alpha
+\bigl(\phi(m,j,k)-un(-2k-j-1+num+2n)\bigr)\delta.
\end{equation*}
Hence
\begin{equation*}
j'\in [0,um]\cap (j+2\mathbb{Z})
\quad\text{and}\quad
k=-\frac{j'+j}{2}-1+n(um+2).
\end{equation*}
In this case,
\begin{equation*}
p(\lambda_k)=-1
\quad\text{and}\quad
\{\lambda_k\}
=m\Lambda_0+\frac{j'\alpha}{2}
+\bigl(\phi(m,j,k)-un(j'+1-num-2n)\bigr)\delta.
\end{equation*}

\item[(c)] (Proof of (2.iii).)
Put $M=um+2$. We consider first the case
$k=\frac{j'-j}{2}-nM$.
Write $k=qm+r$ for some $q\in\mathbb{Z}$, $0\leq r<m$. Then
\begin{align*}
-N_k
&=\phi(m,j,r)-q(k+r+j)-un(j'+1-nM) \\
&=\phi(m,j,r)
-\frac{\left(\frac{j'-j}{2}-nM-r\right)
\left(\frac{j+j'}{2}-nM+r\right)}{m}
-un(j'+1-nM) \\
&=n^2M\left(u-\frac{M}{m}\right)
-n\left(u+uj'-\frac{M}{m}j'\right)
+\frac{j^2-j'^2}{4m}
+\left(\frac{r^2}{m}+\frac{rj}{m}+\phi(m,j,r)\right).
\end{align*}
We have
\begin{equation*}
\frac{r^2}{m}+\frac{rj}{m}+\phi(m,j,r)
=\begin{cases}
\frac{1}{m}r(r+j-m), & \text{if } 0\leq r\leq m-j,\\[2pt]
\frac{1}{m}(r-m)(r+j-m), & \text{if } m-j\leq r<m.
\end{cases}
\end{equation*}
The condition $0\leq r\leq m-j$ can be rewritten as
\begin{equation*}
\frac{j-m}{2}
\leq -nM+\frac{j'}{2}-m\left(q+\frac{1}{2}\right)
\leq \frac{m-j}{2},
\end{equation*}
and $m-j\leq r<m$ can be rewritten as
\begin{equation*}
-\frac{j}{2}
\leq -nM+\frac{j'}{2}-m(q+1)
<\frac{j}{2}.
\end{equation*}
This implies that $\frac{r^2}{m}+\frac{rj}{m}+\phi(m,j,r)$ equals
\begin{equation}\label{Pjj1}
\begin{cases}
\displaystyle
\frac{\left|-nM+\frac{j'}{2}-\frac{m}{2}p\right|^2
-\frac{(m-j)^2}{4}}{m},
& \begin{aligned}
\text{if } \exists\, p\in 2\mathbb{Z}+1 \text{ such that }\\
\left|-nM+\frac{j'}{2}-\frac{m}{2}p\right|
\leq \frac{m-j}{2},
\end{aligned}
\\[10pt]
\displaystyle
\frac{\left|-nM+\frac{j'}{2}-\frac{m}{2}p\right|^2
-\frac{j^2}{4}}{m},
& \begin{aligned}
\text{if } \exists\, p\in 2\mathbb{Z} \text{ such that }\\
\left|-nM+\frac{j'}{2}-\frac{m}{2}p\right|
\leq \frac{j}{2}.
\end{aligned}
\end{cases}
\end{equation}

Let $P_{j,j'}:\mathbb{Z}\to\mathbb{R}$ be the function mapping $-n$ to~\eqref{Pjj1}, and let
$F_{j,j'}:\mathbb{Z}\to\mathbb{R}$ be defined by
\begin{equation*}
F_{j,j'}(t)
=t^2M\left(u-\frac{M}{m}\right)
+t\left(u+uj'-\frac{M}{m}j'\right)
+\frac{j^2-j'^2}{4m}
+P_{j,j'}(t).
\end{equation*}
Then $-N_k=F_{j,j'}(-n)$.

We will show that the maximum of $F_{j,j'}(-n)$ occurs at $n=0$.
Let $P:\mathbb{R}\times[0,m]\times[0,um]\to\mathbb{R}$ be defined by
$P(t,j,j')=P_{j,j'}(t)$, and similarly define
$F(t,j,j')=F_{j,j'}(t)$.
The function $F$ is continuous and piecewise smooth.
Set
\[
\Delta(t,j,j')=F(t+1,j,j')-F(t,j,j').
\]
We will prove that $\Delta$ is nonincreasing in $t$ and that
$\Delta(-1,j,j')>0>\Delta(0,j,j')$.
Indeed,
\begin{equation*}
\Delta(t,j,j')
=2tM\left(u-\frac{M}{m}\right)
+(M+j')\left(u-\frac{M}{m}\right)
+u+P(t+1,j,j')-P(t,j,j').
\end{equation*}
Let $p_1$ and $p_0$ be the integers appearing in
$P_{j,j'}(t+1)$ and $P_{j,j'}(t)$, respectively.
By definition, $\frac{p_1-p_0}{2}\geq u$.
Hence, where they exist,
\begin{align*}
\partial_t\Delta(t,j,j')&=
2M\left(u-\frac{p_1-p_0}{2}\right)\leq 0,\\
\partial_{j'}\Delta(t,j,j')&=
u-\frac{p_1-p_0}{2}\leq 0.
\end{align*}
Thus $\Delta$ is nonincreasing in $t$ and $j'$.
Therefore,
$\Delta(0,j,j')\leq \Delta(0,j,0)$ and
$\Delta(-1,j,um)\leq \Delta(-1,j,j')$.
One checks easily that
$\Delta(0,j,0)<0<\Delta(-1,j,um)$.

This implies that $F(0,j,j')>F(t,j,j')$ for all $t\in\mathbb{Z}$,
$t\neq 0$, i.e.\ $F_{j,j'}(-n)$ attains its maximum at $n=0$.
Hence, in the case $k=\frac{j'-j}{2}-nM$, the minimum of $N_k$
occurs at $n=0$.

For the case $k=-\frac{j+j'}{2}-1+nM$, since
$k=\frac{j'-j}{2}+\left(n-\frac{j'+1}{M}\right)M$, we have
\begin{equation*}
-N_k=F\left(n-\frac{j'+1}{M},j,j'\right).
\end{equation*}
Thus $N_k$ attains its minimum at $n=0$ or $n=1$.
\end{enumerate}
\end{proof}

The next proposition will be used to prove Theorem~\ref{b=b.a11}.

\begin{proposition}\label{sol1}
With $\phi(m,j,k)$ defined as in~\eqref{nk1}, for each
$j\in[0,m]$ and $j'\in[0,um]\cap (j+2\mathbb{Z})$, we have
\begin{equation*}
-\phi\!\left(m,j,\frac{j'-j}{2}\right)
=
\min\!\left(
-\phi\!\left(m,j,-\frac{j'+j}{2}-1\right),
u-\phi\!\left(m,j,-\frac{j'+j}{2}+1\right)
\right)
\end{equation*}
if and only if one of the following conditions holds:
\begin{enumerate}
\item[(A1)] $m>1$ and $j'\leq j-2$;
\item[(B1)] $m>1$, $u$ is even, and $j'\geq um-j+1$;
\item[(C1)] $m>1$, $u$ is odd, and $j'\geq m(u-1)+j+2$.
\end{enumerate}
\end{proposition}

\begin{proof}
We use the fact that
\begin{equation}\label{technk1}
\phi(m,j,-(j+k))=\phi(m,j,k).
\end{equation}
Indeed, if
$\Lambda = m\Lambda_0 +\frac{j\alpha}{2}+b\delta \in P^m_+$
and
$\lambda = \Lambda +k\alpha +\phi(m,j,k)\delta \in \max(\Lambda)$,
then
$s_1(\lambda)=\Lambda -(j+k)\alpha + \phi(m,j,k)\delta\in \max(\Lambda)$.

Using~\eqref{technk1}, we rewrite
\begin{equation*}\label{eq1}
\phi\!\left(m,j,-\frac{j'+j}{2}-1\right)
=\phi\!\left(m,j,\frac{j'-j}{2}\right)
\end{equation*}
as $\phi(m,j,x)=\phi(m,j,x+1)$ with $x=-\frac{j'+j}{2}-1$.
Using~\eqref{nk1}, this holds if and only if~(A1) holds.

Similarly, using~\eqref{technk1}, we rewrite
\begin{equation*}\label{equ1}
-u+\phi\!\left(m,j,-\frac{j'+j}{2}+1\right)
\leq \phi\!\left(m,j,\frac{j'-j}{2}\right)
\end{equation*}
as $\phi(m,j,x+1)-u\leq \phi(m,j,x)$ with $x=-\frac{j'+j}{2}$.
Using~\eqref{nk1}, this holds if and only if~(B1) or~(C1) holds.
\end{proof}

We now turn to the proof of Theorem~\ref{b=b.a11}.

\begin{proof}
The first step is to write explicitly $\ch_{\Lambda}$ in
Proposition~\ref{propchar}. This is done by substituting the values of
$p(\lambda_k)$ and $\{\lambda_k\}$ from Proposition~\ref{prop1} into
formula~\eqref{char.p}. We can rewrite $\ch_{\Lambda}$ as
\begin{equation}\label{rewritechar.1.1}
\sum_{\substack{j'\in [0,um]\\ j'\in j+2\mathbb{Z}}}
\dot{\ch}_{m\Lambda_0+\frac{j'\alpha}{2}}
\left(
\sum_{\substack{n\in\mathbb{Z}\\ k=\frac{j'-j}{2}-n(um+2)}}
q^{N_k}c^\Lambda_{\lambda_k}
-
\sum_{\substack{n\in\mathbb{Z}\\ k=-\frac{j'+j}{2}-1+n(um+2)}}
q^{N_k}c^\Lambda_{\lambda_k}
\right).
\end{equation}
Identity~\eqref{rewritechar.1.1} implies the conditions
$m'=m$ and $j'\in [0,um]\cap (j+2\mathbb{Z})$ required in
Theorem~\ref{b=b.a11}.
The coefficients of $c^\Lambda_{\lambda_k}$ are positive integers,
since $\lambda_k\in \max(\Lambda)$.
Proposition~\ref{prop1} shows that $N_k$ attains its minimum at $n=0$
for the left sum in~\eqref{rewritechar.1.1}, and at $n=0$ or $n=1$ for
the right sum.
The corresponding minima are
\begin{equation}\label{min1.1}
-\phi\!\left(m,j,\frac{j'-j}{2}\right)
\quad\text{and}\quad
\min\!\left(
-\phi\!\left(m,j,-\frac{j'+j}{2}-1\right),
u-\phi\!\left(m,j,-\frac{j'+j}{2}+1\right)
\right).
\end{equation}
By~\eqref{ba11}, \eqref{bulessthanb}, \eqref{rewritechar.1.1},
and~\eqref{min1.1}, we obtain
\begin{equation}\label{inequala1}
-\phi\!\left(m,j,\frac{j'-j}{2}\right)
\leq
\min\!\left(
-\phi\!\left(m,j,-\frac{j'+j}{2}-1\right),
u-\phi\!\left(m,j,-\frac{j'+j}{2}+1\right)
\right).
\end{equation}
Moreover, equality in~\eqref{inequala1} holds if and only if one of
(A1), (B1), or (C1) in Proposition~\ref{sol1} is satisfied.
Thus, for any
$\lambda=m\Lambda_0+\frac{j'\alpha}{2}\in\mathcal{A}_u(\Lambda)$,
the inequality~\eqref{inequala1} is strict.
By~\eqref{rewritechar.1.1}, this implies that
$\lambda+\phi(m,j,\frac{j'-j}{2})\delta\in \max_u(\Lambda)$.
Therefore,
\[
b_{\Lambda,\lambda,u}
=\phi\!\left(m,j,\frac{j'-j}{2}\right)
=b_{\Lambda,\lambda}.
\]
\end{proof}
\subsubsection{Computation of $b_{\Lambda,\lambda}$ for the case $A^{(2)}_2$}\label{bAa.A22}
For type $A^{(2)}_2$, the computation parallels the case $A^{(1)}_1$.
Let $\Lambda_0$ be the $0$-th fundamental weight and let $\alpha=\alpha_1$ be the simple root.
Fix $m\in \mathbb{Z}_{>0}$. Then
\begin{equation*}\label{Pm+.2.2}
\overline{P^m_+}
=
\left\lbrace
m\Lambda_0+\frac{j\alpha}{2}
\,\middle|\,
j\in \left[0,\frac{m}{2}\right]\cap\mathbb{Z}
\right\rbrace.
\end{equation*}
We now give an explicit description of $\max(\Lambda)$ and of the number
$b_{\Lambda,\lambda}$ in this case.

\begin{proposition}\label{maxLamdaandnk2}
Let $\Lambda = m\Lambda_0 +\frac{j\alpha}{2} \in \overline{P^m_+}$. For each $k \in \frac{1}{2}\mathbb{Z}$, let $\phi(m,j,k)$ be the number uniquely determined by $k,m,j$ as follows:
\begin{enumerate}
\item[1.] Write $k=\frac{m}{2} q+r$ for some $q\in \mathbb{Z}$ and $r \in [0,\frac{m}{2})$.
\item[2.] Set
\begin{equation}\label{nk2}
\phi(m,j,k) = -q(k+r+j)+
\begin{cases}
-r & \text{if } r\in [0,\frac{m}{2}-j], \\[2pt]
\frac{m}{2}-j-2r & \text{if } r \in [\frac{m}{2}-j,\frac{m}{2})\cap(\frac{m}{2}+\mathbb{Z}), \\[2pt]
\frac{m-1}{2}-j-2r & \text{if } r \in [\frac{m-1}{2}-j,\frac{m}{2})\cap(\frac{m-1}{2}+\mathbb{Z}).
\end{cases}
\end{equation}
\end{enumerate}
Then we have
\begin{equation*}\label{maxLamda.2.2}
\max(\Lambda) =\left\lbrace \Lambda+k\alpha + \phi(m,j,k)\delta \,\middle\vert\, k\in \frac{1}{2}\mathbb{Z}\right\rbrace.
\end{equation*}
Equivalently, for each $\lambda = m\Lambda_0 +\frac{j'\alpha}{2}$ with $j'\in \mathbb{Z}$, we have
\begin{equation}\label{ba22}
b_{\Lambda,\lambda} = \phi\!\left(m,j,\frac{j'-j}{2}\right).
\end{equation}
\end{proposition}

\begin{proof}
We have $\max(\Lambda)=W(\max(\Lambda)\cap P_+)$. By Lemma \ref{max}, the set $\max(\Lambda)\cap P_+$ contains exactly the elements
\begin{equation*}
\Lambda-m_0\alpha_0,\quad \Lambda-m_1\alpha,\quad \Lambda-\alpha_0-m_2\alpha
\end{equation*}
such that
\begin{equation*}
m_i\in\mathbb{Z}_{\geq 0},\quad
m_0 \leq \frac{m}{2}-j,\quad
m_1 \leq \frac{j}{2},\quad
\frac{j+1}{2}-\frac{m}{4}\leq m_2 \leq \frac{j+1}{2}.
\end{equation*}
Recall that
\[
W =\left\lbrace t_{\frac{n\alpha}{2}},\ t_{\frac{n\alpha}{2}}s_1 \,\middle\vert\, n\in \mathbb{Z}\right\rbrace.
\]
We have
\begin{align}
t_{\frac{n\alpha}{2}}(\Lambda-m_0\alpha_0)
&= \Lambda +\frac{mn+m_0}{2}\alpha
-\frac{(mn+2j+2m_0)n+m_0}{2}\delta, \label{tnb1}\\
t_{\frac{n\alpha}{2}}s_1(\Lambda-m_0\alpha_0)
&=\Lambda +\frac{mn-2j-m_0}{2}\alpha
-\frac{(mn-2j-2m_0)n+m_0}{2}\delta, \label{tnb2}\\
t_{\frac{n\alpha}{2}}(\Lambda-m_1\alpha)
&=\Lambda+\frac{mn-2m_1}{2}\alpha
-\frac{(mn+2j-4m_1)n}{2}\delta, \label{tnb3}\\
t_{\frac{n\alpha}{2}}s_1(\Lambda-m_1\alpha)
&=\Lambda +\frac{mn-2j+2m_1}{2}\alpha
-\frac{(mn-2j+4m_1)n}{2}\delta, \label{tnb4}\\
t_{\frac{n\alpha}{2}}(\Lambda-\alpha_0-m_2\alpha)
&=\Lambda+\frac{mn+1-2m_2}{2}\alpha
-\frac{(mn+2j+2-4m_2)n+1}{2}\delta, \label{tnb5}\\
t_{\frac{n\alpha}{2}}s_1(\Lambda-\alpha_0-m_2\alpha)
&=\Lambda+\frac{mn-1-2j+2m_2}{2}\alpha
-\frac{(mn-2j-2+4m_2)n+1}{2}\delta. \label{tnb6}
\end{align}
Thus, an element $\lambda$ in $\max(\Lambda)$ has the form $\Lambda+r\alpha+n'\delta$ for some $r,n'\in \frac{1}{2}\mathbb{Z}$.
Fix such a $\lambda$. Then, for any $q\in \mathbb{Z}$,
\begin{equation*}
t_{\frac{q\alpha}{2}}(\Lambda+r\alpha+n'\delta)
= \Lambda + \left(\frac{m}{2}q+r\right)\alpha
+ \left(n' - \left(j+2r+\frac{m}{2}q \right)q\right)\delta
\end{equation*}
is still in $\max(\Lambda)$. Set $k=\frac{m}{2} q+r$ for some $q\in\mathbb{Z}$. Then
\begin{equation*}
\Lambda+k\alpha+n''\delta \in \max(\Lambda),
\qquad\text{where } n''=n'-q(k+r+j).
\end{equation*}
Assume now that $0\leq r < \frac{m}{2}$. Then the expression $k=\frac{m}{2} q+r$ is the Euclidean division. By
(\ref{tnb1})--(\ref{tnb6}), we obtain
\begin{equation*}
n'=
\begin{cases}
-r & \text{if } r\in [0,\frac{m}{2}-j], \\[2pt]
\frac{m}{2}-j-2r & \text{if } r\in [\frac{m}{2}-j,\frac{m}{2}) \cap (\frac{m}{2}+\mathbb{Z}), \\[2pt]
\frac{m-1}{2}-j-2r & \text{if } r\in [\frac{m-1}{2}-j,\frac{m}{2})\cap(\frac{m-1}{2}+\mathbb{Z}).
\end{cases}
\end{equation*}
Hence, we obtain $n''= \phi(m,j,k)$ as in (\ref{nk2}). This implies
\begin{equation*}
\max(\Lambda)
=\left\lbrace \Lambda+k\alpha+ \phi(m,j,k)\delta \,\middle\vert\, k\in \frac{1}{2}\mathbb{Z}\right\rbrace.
\end{equation*}
\end{proof}
\subsubsection{The conditions for $b_{\Lambda,\lambda,u}=b_{\Lambda,\lambda}$ in the case $A^{(2)}_2$}

As in the case $A^{(1)}_1$, we give conditions ensuring
$b_{\Lambda,\lambda,u}=b_{\Lambda,\lambda}$.
Together with Theorem~\ref{BAa}, this yields a description of the support
$\Gamma(\mathfrak{g},\mathfrak{g}[u])$ for type $A^{(2)}_2$.

\begin{theorem}\label{b=b.a22}
Let $\mathfrak{g}$ be the affine Kac--Moody algebra of type $A^{(2)}_2$. Let $\Lambda = m\Lambda_0 + \frac{j\alpha}{2} \in P_+$ and let $\lambda = m'\Lambda_0 + \frac{j'\alpha}{2} \in \dot{P}_+$.
\begin{enumerate}
    \item[1.] If there exists $b\in \mathbb{C}$ such that $(\Lambda,\lambda+b\delta)$ belongs to $\Gamma(\mathfrak{g},\mathfrak{g}[u])$, then $m'=m$.
    \item[2.] If moreover
    \begin{enumerate}
        \item[a.] $j \leq j'$; and
        \item[b.] $j' \in \frac{m(u-1)}{2}-j + \bigl(2\mathbb{Z}_{\geq 0} \cup \mathbb{Z}_{<0}\bigr)$,
    \end{enumerate}
    then $b_{\Lambda,\lambda,u} = b_{\Lambda,\lambda}$.
\end{enumerate}
\end{theorem}

For $\Lambda$ as in Theorem~\ref{b=b.a22}, let $\mathcal{A}_u(\Lambda)$ denote the set of
$\lambda$ for which the theorem applies, namely
\begin{equation}\label{AuA.22}
\mathcal{A}_u(\Lambda)
=\left\lbrace
m\Lambda_0+\frac{j'\alpha}{2}
\,\middle\vert\,
j'\in \left[j,\frac{um}{2}\right]\cap \mathbb{Z}
\cap \left(
\frac{m(u-1)}{2}-j+(2\mathbb{Z}_{\geq 0}\cup \mathbb{Z}_{<0})
\right)
\right\rbrace.
\end{equation}

\begin{remark}
This result extends \cite[Theorem~2.2]{KacWaki} for $A^{(2)}_2$ from level one to arbitrary positive level.
\end{remark}

We now state auxiliary propositions.

\begin{proposition}\label{prop2}
Let $\mathfrak{g}$ be the affine Kac--Moody algebra of type $A^{(2)}_2$. Fix $m \in \mathbb{Z}_{>0}$ and $u \in \mathbb{Z}_{>1}$ relatively prime to $2$. Let
$\Lambda = m\Lambda_0 + \frac{j\alpha}{2} \in \overline{P^m_+}$.
\begin{enumerate}
\item[1.] We parametrize $\lambda \in \max(\Lambda)$ such that $\lambda+\dot{\rho}$ is regular with respect to $\dot{W}$ by
\[
\lambda_k = \Lambda + k\alpha + \phi(m,j,k)\delta.
\]
Then the only possible values of $k$ are
\begin{equation*}
k=\frac{j'-j}{2}-n\frac{um+3}{2}
\quad \text{and} \quad
k=-\frac{j'+j}{2}-1+n\frac{um+3}{2},
\end{equation*}
where $j'\in \left[0,\frac{um}{2}\right]\cap \mathbb{Z}$ and $n\in \mathbb{Z}$.
\item[2.] Let
\begin{equation*}
N_k=-\phi(m,j,k)+un\left(j'+1-n\frac{um+3}{2}\right).
\end{equation*}
Then
\begin{enumerate}
\item[i.] If $k=\frac{j'-j}{2}-n\frac{um+3}{2}$, then $p(\lambda_k)=1$ and
\[
\{\lambda_k\}=m\Lambda_0+\frac{j'\alpha}{2}-N_k\delta.
\]
\item[ii.] If $k=-\frac{j'+j}{2}-1+n\frac{um+3}{2}$, then $p(\lambda_k)=-1$ and
\[
\{\lambda_k\}=m\Lambda_0+\frac{j'\alpha}{2}-N_k\delta.
\]
\item[iii.] The function $N_k$, considered as a function of $n$, attains its minimum at $n=0$ in the first case and at $n=0$ or $n=1$ in the second case.
\end{enumerate}
\end{enumerate}
\end{proposition}

\begin{proof}
We need the following data:
\begin{equation*}
|\alpha|^2=4,\quad 
\dot{\rho} = \frac{3}{u}\Lambda_0 + \frac{1}{2}\alpha,\quad 
\dot{W}=\left\lbrace t_{\frac{un\alpha}{2}},\, t_{\frac{un\alpha}{2}}s_1 \,\middle\vert\, n\in \mathbb{Z}\right\rbrace.
\end{equation*}
Since $\lambda_k+\dot{\rho}$ is regular with respect to $\dot{W}$, there exists a unique $\sigma \in \dot{W}$ and
\[
\mu = m\Lambda_0 + \frac{j'\alpha}{2} + b'\delta \in \dot{P}^{um}_+
\]
such that $\sigma(\lambda_k+\dot{\rho})=\mu + \dot{\rho}$.

\begin{enumerate}
\item[a.] (Proof of 1. and 2.i.) If $\sigma = t_{\frac{un\alpha}{2}}$ for some $n\in \mathbb{Z}$, then
$\sigma(\lambda_k+\dot{\rho})-\dot{\rho}$ equals
\begin{equation*}
m\Lambda_0 + \left(k+\frac{j}{2}+n\frac{um+3}{2}\right)\alpha
+ \left(\phi(m,j,k)-un \left(2k+j+1+n\frac{um+3}{2}\right)\right)\delta.
\end{equation*}
Hence,
\begin{equation*}
j' \in \left[0,\frac{um}{2}\right] \cap \mathbb{Z}
\quad \text{and} \quad
k=\frac{j'-j}{2}-n\frac{um+3}{2}.
\end{equation*}
In this case, we have
\begin{equation*}
p(\lambda_k)=1
\quad \text{and} \quad
\{\lambda_k\}
= m\Lambda_0 + \frac{j'\alpha}{2}
+ \left(\phi(m,j,k)-un \left(j'+1-n\frac{um+3}{2}\right)\right)\delta.
\end{equation*}

\item[b.] (Proof of 1. and 2.ii.)
If $\sigma = t_{\frac{un\alpha}{2}}s_1$ for some $n\in \mathbb{Z}$, then
$\sigma(\lambda_k+\dot{\rho})-\dot{\rho}$ equals
\begin{equation*}
m\Lambda_0
+ \left(-k-\frac{j}{2}-1+n\frac{um+3}{2}\right)\alpha
+ \left(\phi(m,j,k)-un \left(-2k-j-1+n\frac{um+3}{2}\right)\right)\delta.
\end{equation*}
Hence,
\begin{equation*}
j' \in \left[0,\frac{um}{2}\right] \cap \mathbb{Z}
\quad \text{and} \quad
k=-\frac{j'+j}{2}-1+n\frac{um+3}{2}.
\end{equation*}
In this case, we have
\begin{equation*}
p(\lambda_k)=-1
\quad \text{and} \quad
\{\lambda_k\}
= m\Lambda_0 + \frac{j'\alpha}{2}
+ \left(\phi(m,j,k)-un \left(j'+1-n\frac{um+3}{2}\right)\right)\delta.
\end{equation*}
\item[c.] (Proof of 2.iii.)
Put $M=um+3$. We consider the first case
\[
k=\frac{j'-j}{2}-\frac{nM}{2}.
\]
Write $k=\frac{m}{2}q+r$ for some $q\in\mathbb{Z}$, $0\leq r<\frac{m}{2}$. Then
\begin{align*}
-N_k
&= \phi(m,j,r)-q(k+r+j)-un\left(j'+1-\frac{nM}{2}\right) \\
&=\phi(m,j,r)-\frac{2\left(\frac{j'-j}{2}-\frac{nM}{2}-r\right)
\left(\frac{j+j'}{2}-\frac{nM}{2}+r\right)}{m}
-un\left(j'+1-\frac{nM}{2}\right)\\
&=n^2\frac{M}{2}\left(u-\frac{M}{m}\right)
-n\left(u+uj'-\frac{M}{m}j'\right)
+\frac{j^2-j'^2}{2m}
+\left(\frac{2r^2}{m}+\frac{2rj}{m}+\phi(m,j,r)\right).
\end{align*}

We have $\frac{2r^2}{m}+\frac{2rj}{m}+\phi(m,j,r)$ equal to
\begin{equation*}
\begin{cases}
\frac{2}{m} r\left(r+j-\frac{m}{2}\right)
& \text{if } 0\leq r\leq \frac{m}{2}-j, \\[2pt]
\frac{2}{m}\left(r-\frac{m}{2}\right)\left(r+j-\frac{m}{2}\right)
& \text{if } \frac{m}{2}-j\leq r<\frac{m}{2},\ r\in \frac{m}{2}+\mathbb{Z},\\[2pt]
\frac{2}{m}\left(r-\frac{m}{2}\right)\left(r+j-\frac{m}{2}\right)-\frac{1}{2}
& \text{if } \frac{m-1}{2}-j\leq r<\frac{m}{2},\ r\in \frac{m+1}{2}+\mathbb{Z}.
\end{cases}
\end{equation*}

The condition $0\leq r\leq \frac{m}{2}-j$ can be rewritten as
\begin{equation*}
\frac{2j-m}{4}
\leq \frac{-nM+j'}{2}-\frac{m}{2}\left(q+\frac{1}{2}\right)
\leq \frac{m-2j}{4},
\end{equation*}
and $\frac{m}{2}-j\leq r<\frac{m}{2}$ can be rewritten as
\begin{equation*}
\frac{-j}{2}
\leq \frac{-nM+j'}{2}-\frac{m}{2}(q+1)
< \frac{j}{2}.
\end{equation*}

It implies that $\frac{2r^2}{m}+\frac{2rj}{m}+\phi(m,j,r)$ equals
\begin{equation}\label{Pjj2}
\begin{cases}
\frac{2}{m}\left(
\left|\frac{-nM+j'}{2}-\frac{m}{4}p\right|^2
-\left|\frac{m-2j}{4}\right|^2
\right)
& \text{if $\exists p\in 2\mathbb{Z}+1$ such that} \\[-2pt]
& \left|\frac{-nM+j'}{2}-\frac{m}{4}p\right|\leq \frac{m-2j}{4}, \\[6pt]
\frac{2}{m}\left(
\left|\frac{-nM+j'}{2}-\frac{m}{4}p\right|^2
-\frac{j^2}{4}
\right)
& \text{if $\exists p\in 2\mathbb{Z}$ such that} \\[-2pt]
& \left|\frac{-nM+j'}{2}-\frac{m}{4}p\right|\leq \frac{j}{2}, \\[-2pt]
& \text{$\frac{j'-j}{2}-\frac{nM}{2}-\frac{m}{2}\left(\frac{p}{2}-1\right)
\in \frac{m}{2}+\mathbb{Z}$}, \\[6pt]
\frac{2}{m}\left(
\left|\frac{-nM+j'}{2}-\frac{m}{4}p\right|^2
-\frac{j^2}{4}
\right)-\frac{1}{2}
& \text{if $\exists p\in 2\mathbb{Z}$ such that} \\[-2pt]
& \frac{-nM+j'}{2}-\frac{m}{4}p
\in \left[-\frac{1+j}{2},\frac{j}{2}\right], \\[-2pt]
& \text{$\frac{j'-j}{2}-\frac{nM}{2}-\frac{m}{2}\left(\frac{p}{2}-1\right)
\in \frac{m+1}{2}+\mathbb{Z}$}.
\end{cases}
\end{equation}

Let $P_{j,j'}:\mathbb{Z}\rightarrow \mathbb{R}$ be the function mapping $-n$ to
(\ref{Pjj2}). Let $F_{j,j'}:\mathbb{Z}\rightarrow \mathbb{R}$ be defined by
\begin{equation*}
F_{j,j'}(t)
= t^2\frac{M}{2}\left(u-\frac{M}{m}\right)
+ t\left(u+uj'-\frac{M}{m}j'\right)
+ \frac{j^2-j'^2}{2m}
+ P_{j,j'}(t).
\end{equation*}
Thus, $-N_k=F_{j,j'}(-n)$.

We will show that the maximum of $F_{j,j'}(n)$ occurs when $n=0$, i.e.,
$k=\frac{j'-j}{2}$. To do this, we show that the upper bound function
$F^+:\mathbb{R}\times[0,\frac{m}{2}]\times[0,\frac{um}{2}]\rightarrow\mathbb{R}$
and the lower bound function
$F^-:\mathbb{R}\times[0,\frac{m}{2}]\times[0,\frac{um}{2}]\rightarrow\mathbb{R}$
of $F_{j,j'}(t)$, defined below, attain their maxima over $t\in\mathbb{Z}$ at
$t=0$.

\begin{align*}
F^+(t,j,j')
&= t^2\frac{M}{2}\left(u-\frac{M}{m}\right)
+ t\left(u+uj'-\frac{M}{m}j'\right)
+ \frac{j^2-j'^2}{2m}
+ P^+(t,j,j'),\\
F^-(t,j,j')
&= t^2\frac{M}{2}\left(u-\frac{M}{m}\right)
+ t\left(u+uj'-\frac{M}{m}j'\right)
+ \frac{j^2-j'^2}{2m}
+ P^-(t,j,j').
\end{align*}

We consider only $F^+$; the argument for $F^-$ is similar.
The function $F^+$ is continuous and piecewise smooth. Let
\[
\Delta^+(t,j,j')=F^+(t+1,j,j')-F^+(t,j,j').
\]
Then $\Delta^+(t,j,j')$ is nonincreasing in $t$, and
$\Delta^+(-1,j,j')>0>\Delta^+(0,j,j')$. Indeed,
\begin{equation*}
\Delta^+(t,j,j')
= tM\left(u-\frac{M}{m}\right)
+ \left(j'+\frac{M}{2}\right)\left(u-\frac{M}{m}\right)
+ u
+ P^+(t+1,j,j')-P^+(t,j,j').
\end{equation*}

Denote by $p_1$ and $p_0$ the integers defining $P^+(t+1,j,j')$ and
$P^+(t,j,j')$, respectively. By definition,
$\frac{p_1-p_0}{2}\geq u$. Hence, where they exist,
\begin{align*}
\partial_t \Delta^+(t,j,j')
&= M\left(u-\frac{p_1-p_0}{2}\right)\leq 0,\\
\partial_{j'} \Delta^+(t,j,j')
&= u-\frac{p_1-p_0}{2}\leq 0.
\end{align*}

Thus, $\Delta^+$ is nonincreasing in $t$ and $j'$. Therefore,
$\Delta^+(0,j,j')\leq \Delta^+(0,j,0)$ and
$\Delta^+(-1,j,\frac{um}{2})\leq \Delta^+(-1,j,j')$.
A direct check shows that
$\Delta^+(0,j,0)<0<\Delta^+(-1,j,\frac{um}{2})$.

It follows that $F^+(0,j,j')>F^+(t,j,j')$ for all $t\in\mathbb{Z}$,
$t\neq 0$. The same result holds for $F^-$. Hence,
$F_{j,j'}(-n)$ attains its maximum when $n=0$. Consequently, in the case
$k=\frac{j'-j}{2}-\frac{nM}{2}$, the minimum of $N_k$ occurs when $n=0$.

For the case
\[
k=-\frac{j'+j}{2}-1+\frac{nM}{2},
\]
since
\[
k=\frac{j'-j}{2}+\left(n-2\frac{j'+1}{M}\right)\frac{M}{2},
\]
we have
\begin{equation*}
-N_k=F\left(n-2\frac{j'+1}{M},j,j'\right).
\end{equation*}
Thus, $N_k$ attains its minimum when $n=0$ or $1$.
\end{enumerate}
\end{proof}

The next proposition will be used to prove Theorem \ref{b=b.a22}.

\begin{proposition}\label{sol2}
With $\phi(m,j,k)$ defined as in (\ref{nk2}), for each
$j \in \left[0,\frac{m}{2}\right]$ and
$j' \in \left[0,\frac{um}{2}\right]\cap \mathbb{Z}$, we have
\begin{equation*}
-\phi\!\left(m,j,\frac{j'-j}{2}\right)
=
\min\!\left(
-\phi\!\left(m,j,-\frac{j'+j}{2}-1\right),
\frac{u}{2}-\phi\!\left(m,j,-\frac{j'+j}{2}+\frac{1}{2}\right)
\right)
\end{equation*}
if and only if one of the following two conditions is satisfied:
\begin{enumerate}
\item[(A2)]
$m>2$ and $j'\leq j-1$.
\item[(B2)]
$m\geq 2$ and
$j'\in \frac{m(u-1)}{2}+1-j+2\mathbb{Z}_{\geq 0}$.
\end{enumerate}
\end{proposition}

\begin{proof}
We again use the equality (\ref{technk1}) to rewrite
\begin{equation*}
\phi\!\left(m,j,-\frac{j'+j}{2}-1\right)
=
\phi\!\left(m,j,\frac{j'-j}{2}\right),
\qquad
-\frac{u}{2}
+\phi\!\left(m,j,-\frac{j'+j}{2}+\frac{1}{2}\right)
=
\phi\!\left(m,j,\frac{j'-j}{2}\right)
\end{equation*}
as
$\phi(m,j,x)=\phi(m,j,x-1)$ and
$\phi(m,j,x)=-\frac{u}{2}+\phi\!\left(m,j,x+\frac{1}{2}\right)$,
where $x=-\frac{j'+j}{2}$.
Using (\ref{nk2}), we can check that this occurs if and only if
condition \textup{(A2)} or \textup{(B2)} is satisfied.
\end{proof}

We now turn to the proof of Theorem \ref{b=b.a22}.

\begin{proof}
The strategy is the same as in the case $A^{(1)}_1$.
The first step is to write $\ch_{\Lambda}$ explicitly using
Proposition \ref{propchar}.
By substituting the values of $p(\lambda_k)$ and $\{\lambda_k\}$
from Proposition \ref{prop2} into the formula (\ref{char.p}),
we can rewrite $\ch_{\Lambda}$ as
\begin{equation}\label{rewritechar.2.2}
\sum\limits_{j'\in \left[0,\frac{um}{2}\right]\cap \mathbb{Z}}
\dot{\ch}_{m\Lambda_0+\frac{j'\alpha}{2}}
\left(
\sum\limits_{\substack{n\in\mathbb{Z}\\
k=\frac{j'-j}{2}-n\frac{um+3}{2}}}
q^{N_k}c^\Lambda_{\lambda_k}
-
\sum\limits_{\substack{n\in\mathbb{Z}\\
k=-\frac{j'+j}{2}-1+n\frac{um+3}{2}}}
q^{N_k}c^\Lambda_{\lambda_k}
\right).
\end{equation}

The identity (\ref{rewritechar.2.2}) implies the condition
$j'\in \left[0,\frac{um}{2}\right]\cap \mathbb{Z}$ required in
Theorem \ref{b=b.a22}.
The coefficients of $c^\Lambda_{\lambda_k}$ in
(\ref{rewritechar.2.2}) are always positive integers, since
$\lambda_k\in \max(\Lambda)$.

Proposition \ref{prop2} shows that the number $N_k$ attains its
minimum at $n=0$ for the terms on the left-hand side of
(\ref{rewritechar.2.2}), and at $n=0$ or $n=1$ for the terms on the
right-hand side.
The corresponding minima of $N_k$ are
\begin{equation}\label{min2.2}
-\phi\!\left(m,j,\frac{j'-j}{2}\right)
\quad\text{and}\quad
\min\!\left(
-\phi\!\left(m,j,-\frac{j'+j}{2}-1\right),
\frac{u}{2}-\phi\!\left(m,j,-\frac{j'+j}{2}+\frac{1}{2}\right)
\right).
\end{equation}

By (\ref{ba22}), (\ref{bulessthanb}), (\ref{rewritechar.2.2}),
and (\ref{min2.2}), we obtain
\begin{equation}\label{inequala2}
-\phi\!\left(m,j,\frac{j'-j}{2}\right)
\leq
\min\!\left(
-\phi\!\left(m,j,-\frac{j'+j}{2}-1\right),
\frac{u}{2}-\phi\!\left(m,j,-\frac{j'+j}{2}+\frac{1}{2}\right)
\right).
\end{equation}

Moreover, equality holds if and only if condition \textup{(A2)} or
\textup{(B2)} in Proposition \ref{sol2} is satisfied.
Thus, for any
$\lambda=m\Lambda_0+\frac{j'\alpha}{2}\in \mathcal{A}_u(\Lambda)$,
the inequality in (\ref{inequala2}) is strict.

By (\ref{rewritechar.2.2}), in this case we have
$\lambda+\phi\!\left(m,j,\frac{j'-j}{2}\right)\delta
\in \max_u(\Lambda)$.
This means that
\[
b_{\Lambda,\lambda,u}
=
\phi\!\left(m,j,\frac{j'-j}{2}\right)
=
b_{\Lambda,\lambda}.
\]
\end{proof}
\subsubsection{Relation between $\Gamma(\mathfrak{g},\mathfrak{g}[u])$ and its saturated setting}

The saturated setting of the support $\Gamma(\mathfrak{g},\mathfrak{g}[u])$ is defined by
\begin{equation*}
\tilde{\Gamma}(\mathfrak{g},\mathfrak{g}[u])
=
\bigl\{
(\Lambda,\lambda) \in P_+ \times \dot{P}_+
\;\big|\;
\lambda \in \Lambda+Q,\ 
\dot{L}(N\lambda) \subset L(N\Lambda)
\text{ for some integer } N>1
\bigr\}.
\end{equation*}

\begin{corollary}\label{satured}
Let $\mathfrak{g}$ be an affine Kac--Moody algebra of type $A^{(1)}_1$ or $A^{(2)}_2$. Fix $u \in \mathbb{Z}_{>1}$ ($u$ is an odd integer in the case of $A^{(2)}_2$). Let $\Lambda \in P_+$ and let $\lambda \in \mathcal{A}_u(\Lambda) \cap (\Lambda+Q)$. For all $b \in \mathbb{C}$, we have
\begin{enumerate}
    \item[1.] $(\Lambda,\lambda+b\delta) \in \tilde{\Gamma}(\mathfrak{g},\mathfrak{g}[u]) \Leftrightarrow d(\Lambda,\lambda+b\delta) \in \Gamma(\mathfrak{g},\mathfrak{g}[u])$ for all $d \in \mathbb{Z}_{\geq 2}$.
    \item[2.] If, in addition, ${h}_{\Lambda,\lambda}^{[u]} \neq 0$, then
    \[
    (\Lambda,\lambda+b\delta) \in \tilde{\Gamma}(\mathfrak{g},\mathfrak{g}[u]) \Leftrightarrow (\Lambda,\lambda+b\delta) \in \Gamma(\mathfrak{g},\mathfrak{g}[u]).
    \]
\end{enumerate}
\end{corollary}

Before the proof, we recall \cite[Lemmas~6.3 and~8.5]{KumarBrown}.

\begin{lemma}\label{lemm6385}
Let $\mathfrak{g}$ be an affine Kac--Moody algebra of type
$A^{(1)}_1$ or $A^{(2)}_2$.
Let $\Lambda \in P_+$ and $\lambda \in \Lambda + Q$.
Fix a positive integer $N$.
Then $\lambda \in \max(\Lambda)$ if and only if
$N\lambda \in \max(N\Lambda)$.
\end{lemma}

We now prove Corollary~\ref{satured}.

\begin{proof}
Fix $d\in \mathbb{Z}_{>0}$. Since $\lambda\in \mathcal{A}_u(\Lambda)\cap(\Lambda+Q)$, we have
$d\lambda\in \mathcal{A}_u(d\Lambda)\cap(d\Lambda+Q)$. By
Theorems~\ref{b=b.a11}, \ref{b=b.a22}, and Lemma~\ref{lemm6385},
\begin{equation}\label{db=db}
b_{d\Lambda,d\lambda,u}=b_{d\Lambda,d\lambda}=d\,b_{\Lambda,\lambda}.
\end{equation}

Assume $(\Lambda,\lambda+b\delta)\in \tilde{\Gamma}(\mathfrak{g},\mathfrak{g}[u])$.
Then $\lambda+b\delta\in \Lambda+Q$ and
$\dot{L}(N\lambda+Nb\delta)\subset L(N\Lambda)$ for some $N>1$.
Hence $b\in \mathbb{Z}$ and, by \eqref{db=db},
\[
Nb\le b_{N\Lambda,N\lambda,u}=N b_{\Lambda,\lambda},
\]
so $b\le b_{\Lambda,\lambda}$.

Therefore, for any $d\ge 2$,
\[
b_{d\Lambda,d\lambda}-db=d(b_{\Lambda,\lambda}-b)\in \mathbb{Z}_{\ge 0}\setminus\{1\}.
\]
By Theorem~\ref{BAa}(1), this yields
$d(\Lambda,\lambda+b\delta)\in \Gamma(\mathfrak{g},\mathfrak{g}[u])$.
If $h_{\Lambda,\lambda}^{[u]}\neq 0$, then by Theorem~\ref{BAa}(2),
$(\Lambda,\lambda+b\delta)\in \Gamma(\mathfrak{g},\mathfrak{g}[u])$.
\end{proof}

\begin{remark}
For $A^{(1)}_1$, we have
\[
\mathcal{A}_u(\Lambda)\cap(\Lambda+Q)=\mathcal{A}_u(\Lambda).
\]
For $A^{(2)}_2$, the set $\mathcal{A}_u(\Lambda)\cap(\Lambda+Q)$ consists of those
elements of $\mathcal{A}_u(\Lambda)$ with the additional constraint
$j'\in j+2\mathbb{Z}$, cf.~\eqref{AuA.22}.
\end{remark}


\section*{Acknowledgments}
\addcontentsline{toc}{section}{Acknowledgments}

The author gratefully thanks his supervisors, Prof.\ Nicolas Ressayre and Prof.\ Kenji Iohara, for proposing the topic, for numerous helpful comments, and for their careful reading of the manuscript. He also thanks Prof.\ Shrawan Kumar, Prof.\ Peter Littelmann, and Prof.\ St\'{e}phane Gaussent for valuable discussions. Further thanks go to Prof.\ C\'{e}dric Lecouvey and Prof.\ Shrawan Kumar, the referees of the thesis, for their thorough reading and for pointing out several important errors, which led to substantial improvements.


\appendix
\section{Coset construction}\label{app}
In this section, we present some facts about the Virasoro algebras.

\subsection{Sugawara construction of Virasoro operators}

We recall the Sugawara construction of the Virasoro operators in \cite{Kac, KacWaki2}.

\subsubsection{Untwisted case}

We recall the Sugawara construction for untwisted affine Kac--Moody algebras.
Let $\mathfrak{g}$ be of type $X_N^{(1)}$ and let $\overline{\mathfrak{g}}$ be the
simple Lie algebra of type $X_N$. Then
\[
\mathfrak{g}
=
\mathbb{C}[t,t^{-1}] \otimes \overline{\mathfrak{g}}
\oplus \mathbb{C}K \oplus \mathbb{C}d,
\]
with Lie bracket
\[
[t^i\!\otimes x+\lambda K+\mu d,\; t^j\!\otimes y+\lambda' K+\mu' d]
=
t^{i+j}\!\otimes[x,y]_0
+\mu j\, t^j\!\otimes y
-\mu' i\, t^i\!\otimes x
+i\delta_{i+j,0}(x|y)_0 K,
\]
for $i,j\in\mathbb{Z}$, $x,y\in\overline{\mathfrak{g}}$, and
$\lambda,\lambda',\mu,\mu'\in\mathbb{C}$.
Write $x^{(m)}=t^m\otimes x$.
Let $\{u_i\}$ and $\{u^i\}$ be dual bases of $\overline{\mathfrak{g}}$,
$(u_i|u^j)=\delta_{ij}$.
The Sugawara operators $T_n$ $(n\in\mathbb{Z})$ are defined by
\[
T_0=\sum_i u_i u^i + 2\sum_{m>0}\sum_i u_i^{(-m)}u^{i(m)},
\qquad
T_n=\sum_{m\in\mathbb{Z}}\sum_i u_i^{(-m)}u^{i(m+n)},\; n\neq0.
\]

The sums are infinite and do not define elements of $\mathcal{U}(\mathfrak{g})$.
They act, however, as well-defined endomorphisms on any
$\mathfrak{g}$-module $V$ in category $\mathcal{O}$, since for each $v\in V$
only finitely many terms act nontrivially.

Let $V$ be a $\mathfrak{g}'$-module in category $\mathcal{O}$ such that $K$
acts by the scalar $k$ with $k+\mathrm{h}^\vee\neq0$.
As endomorphisms of $V$, the following result holds.

\begin{proposition}\label{Tn1}
For any $n,m \in \mathbb{Z}$ and $x \in \overline{\mathfrak{g}}$, we have
\begin{enumerate}
\item[1.] $[x^{(m)}, T_n] = 2(k+\mathrm{h}^\vee) m x^{(m+n)}$.
\item[2.] $[T_m, T_n] = 2(k+\mathrm{h}^\vee)\bigl((m-n)T_{m+n}
+ \delta_{m+n,0}\tfrac{m^3-m}{6}(\dim \overline{\mathfrak{g}})k\,\mathrm{Id}_V\bigr)$.
\end{enumerate}
\end{proposition}

Let $\mathfrak{g}'=[\mathfrak{g},\mathfrak{g}]
=\mathbb{C}[t,t^{-1}]\otimes_{\mathbb{C}}\overline{\mathfrak{g}}\oplus\mathbb{C}K$.
We define the \black{Virasoro operators} on $V$ by
\begin{equation*}
L_n=\frac{1}{2(k+\mathrm{h}^\vee)}\,T_n,\qquad n\in\mathbb{Z},
\end{equation*}
and let $Z$ act by the scalar
\[
\frac{k\,\dim\overline{\mathfrak{g}}}{k+\mathrm{h}^\vee}.
\]
The following proposition is an immediate consequence of Proposition~\ref{Tn1}.

\begin{proposition}\label{Ln1}
For any $n,m \in \mathbb{Z}$ and $x \in \overline{\mathfrak{g}}$, we have
\begin{enumerate}
\item[1.] $[x^{(m)}, L_n] = m x^{(m+n)}$.
\item[2.] $[L_m, L_n] = (m-n)L_{m+n}
+ \delta_{m+n,0}\tfrac{m^3-m}{12} Z$.
\end{enumerate}
\end{proposition}
\subsubsection{Twisted case}

In this part, we recall the Sugawara construction for twisted affine Kac--Moody algebras.
Let $\mathfrak{g}$ be of twisted type $X_N^{(r)}$ with $r=2,3$.
Let $\overline{\mathfrak{g}}$ be the simple Lie algebra of type $X_N$, and let $\sigma$ be its diagram automorphism.
The induced $\mathbb{Z}/r\mathbb{Z}$-grading is
\[
\overline{\mathfrak{g}}=\bigoplus_{\overline{j}\in\mathbb{Z}/r\mathbb{Z}}\overline{\mathfrak{g}}_{\overline{j}},
\qquad
\overline{\mathfrak{g}}_{\overline{j}}
=\{x\in\overline{\mathfrak{g}}\mid \sigma(x)=e^{2\pi ij/r}x\}.
\]
Then
\begin{equation*}
\mathfrak{g}
=\bigoplus_{j\in\mathbb{Z}} t^j\otimes\overline{\mathfrak{g}}_{\overline{j}}
\oplus\mathbb{C}K\oplus\mathbb{C}d',
\end{equation*}
with Lie bracket
\[
[t^i\otimes x+\lambda K+\mu d',\, t^j\otimes y+\lambda'K+\mu'd']
= t^{i+j}\otimes[x,y]_0
+\mu j\,t^j\otimes y-\mu' i\,t^i\otimes x
+\frac{i}{r}\delta_{i+j,0}(x|y)_0 K.
\]
The normalized invariant form on $\mathfrak{g}$ is given by
\[
(t^i\otimes x+\lambda K+\mu d'\mid t^j\otimes y+\lambda'K+\mu'd')'
=\frac{1}{r}\delta_{i+j,0}(x|y)_0+(\lambda'\mu+\lambda\mu').
\]

\begin{remark}
\begin{enumerate}
    \item[1.] We have $(\alpha|\alpha)_0 = 2$ for all long roots $\alpha$ of $\overline{\mathfrak{g}}$. However, $(\alpha|\alpha)' = 2r$ for all long roots $\alpha$ of $\mathfrak{g}$.
    \item[2.] The element $d'$ is not the scaling element of $\mathfrak{g}$ as in the untwisted case. The formula for the scaling element involves $d'$.
\end{enumerate}
\end{remark}

Let $\{u_{i,-\overline{j}}\}$ be a basis of $\overline{\mathfrak{g}}_{-\overline{j}}$ and
$\{u^{i,\overline{j}}\}$ the dual basis of $\overline{\mathfrak{g}}_{\overline{j}}$.
The \black{Sugawara operators} $T_n$ ($n\in\mathbb{Z}$) are defined by
\begin{equation*}
T_0=
\sum_i u_{i,\overline{0}}u^{i,\overline{0}}
+2\sum_{m>0}\sum_i u_{i,-\overline{m}}^{(-m)}u^{i,\overline{m}(m)}
+\frac{r-1}{2r}\bigl(\dim\overline{\mathfrak{g}}-\dim\overline{\mathfrak{g}}_{\overline{0}}\bigr)K,
\end{equation*}
and
\begin{equation*}
T_n=\sum_{m\in\mathbb{Z}}\sum_i
u_{i,-\overline{m}}^{(-m)}u^{i,\overline{m}(m+rn)},
\qquad n\neq0.
\end{equation*}
For any $\mathfrak{g}$-module $V$ in the category $\mathcal{O}$, these operators define
well-defined endomorphisms, since for each $v\in V$ only finitely many terms act nontrivially.

Set
\[
\mathfrak{g}'=[\mathfrak{g},\mathfrak{g}]
=\bigoplus_{j\in\mathbb{Z}} t^j\otimes\overline{\mathfrak{g}}_{\overline{j}}\oplus\mathbb{C}K.
\]
Let $V$ be a $\mathfrak{g}'$-module in $\mathcal{O}$ on which $K$ acts by the scalar $k$,
with $k+\mathrm{h}^\vee\neq0$. As endomorphisms of $V$, we obtain the following result.

\begin{proposition}\label{Tnr}
For any $n,m \in \mathbb{Z}$ and $x \in \overline{\mathfrak{g}}_{\overline{m}}$, we have
\begin{enumerate}
\item[1.] $[x^{(m)}, T_n] = \dfrac{2(k+\mathrm{h}^\vee)}{r} m x^{(m+rn)}$.
\item[2.] $[T_m, T_n] = 2(k+\mathrm{h}^\vee)\Bigl((m-n)T_{m+n}
+ \delta_{m+n,0}\dfrac{m^3-m}{6} r (\dim \overline{\mathfrak{g}}) k\,\mathrm{Id}_V\Bigr)$.
\end{enumerate}
\end{proposition}

We define the \black{Virasoro operators} on $V$ by
\begin{equation*}
L_n=\frac{1}{2(k+\mathrm{h}^\vee)}\,T_n,
\qquad n\in\mathbb{Z},
\end{equation*}
and let $Z$ act by the scalar
\[
\frac{rk\,\dim\overline{\mathfrak{g}}}{k+\mathrm{h}^\vee}.
\]
The following proposition follows directly from Proposition~\ref{Tnr}.

\begin{proposition}\label{Ln1r}
For any $n,m \in \mathbb{Z}$ and $x \in \overline{\mathfrak{g}}_{\overline{m}}$, we have
\begin{enumerate}
    \item[1.] $[x^{(m)}, L_n] = \dfrac{m}{r} x^{(m+rn)}$.
    \item[2.] $[L_m, L_n] = (m-n)L_{m+n}
    + \delta_{m+n,0}\dfrac{m^3-m}{12} Z$.
\end{enumerate}
\end{proposition}
\subsubsection{Anti-involution $\omega_0$ on Virasoro operators}

The first parts of Propositions~\ref{Ln1} and~\ref{Ln1r} define a
$\mathfrak{g}'\rtimes \Vir$-module structure on $V$.
Here $\Vir$ denotes the Virasoro algebra generated by $L_n$ ($n\in\mathbb{Z}$)
and $Z$.
We recall the anti-linear anti-involution $\omega_0$ on
$\mathcal{U}(\mathfrak{g})$ from Subsection~\ref{w0}.
Up to a suitable completion of $\mathcal{U}(\mathfrak{g})$ allowing infinite
series, we obtain the following statement.

\begin{proposition}\label{w0LZ}
We have $\omega_0(L_n)=L_{-n}$ for all $n \in \mathbb{Z}$, and $\omega_0(Z)=Z$.
\end{proposition}

\begin{remark}
Proposition~\ref{w0LZ} implies
\begin{equation*}
    \langle L_n v, w \rangle_\Lambda
    = \langle v, \omega_0(L_n) w \rangle_\Lambda
    = \langle v, L_{-n} w \rangle_\Lambda,
\end{equation*}
and
\begin{equation*}
    \langle Z v, w \rangle_\Lambda
    = \langle v, \omega_0(Z) w \rangle_\Lambda
    = \langle v, Z w \rangle_\Lambda.
\end{equation*}
\end{remark}
\subsection{Coset construction of Virasoro operators for winding subalgebras}\label{coset construction}

Let $\mathfrak{g}$ be an affine Kac--Moody algebra of type $X_N^{(r)}$, and let
$\overline{\mathfrak{g}}$ be the simple Lie algebra of type $X_N$.
We recall the coset construction of Virasoro operators for the winding
subalgebras $\mathfrak{g}[u]$, introduced in \cite{KacWaki}, and briefly review
the necessary details.

\subsubsection{Untwisted case}

In the case $r=1$, we have
\[
\mathfrak{g}
=\mathbb{C}[t,t^{-1}]\otimes \overline{\mathfrak{g}}
\oplus \mathbb{C}K \oplus \mathbb{C}d,
\qquad
\mathfrak{g}[u]
=\mathbb{C}[t^u,t^{-u}]\otimes \overline{\mathfrak{g}}
\oplus \mathbb{C}K \oplus \mathbb{C}d.
\]
For $\Lambda\in P^k_+$, the $\mathfrak{g}$-module $L(\Lambda)$ is viewed as a
$\mathfrak{g}[u]$-module of level $uk$.
Let $\{u_i\}$ and $\{u^i\}$ be dual bases of $\overline{\mathfrak{g}}$, and set
$x^{(n)}=t^n\otimes x$ for $x\in\overline{\mathfrak{g}}$, $n\in\mathbb{Z}$.
The \black{Virasoro operators} $L_n$ ($n\in\mathbb{Z}$) and $Z$ on $L(\Lambda)$,
given by the Sugawara construction, are
\[
L_0=\frac{1}{2(k+\mathrm{h}^\vee)}
\Bigl(\sum_i u_i u^i
+2\sum_{m>0}\sum_i u_i^{(-m)}u^{i(m)}\Bigr),
\]
\[
L_n=\frac{1}{2(k+\mathrm{h}^\vee)}
\Bigl(\sum_{m\in\mathbb{Z}}\sum_i
u_i^{(-m)}u^{i(m+n)}\Bigr),
\quad n\neq0,
\]
and $Z$ acts by the scalar
\[
c_k=\frac{k\,\dim(\overline{\mathfrak{g}})}{k+\mathrm{h}^\vee}.
\]

The \black{Virasoro operators} $\dot{L}_n$ ($n\in\mathbb{Z}$) and $\dot{Z}$ on
$\mathfrak{g}[u]$-modules are defined by
\[
\dot{L}_0=\frac{1}{2(uk+\mathrm{h}^\vee)}
\Bigl(\sum_i u_i u^i
+2\sum_{m>0}\sum_i u_i^{(-um)}u^{i(um)}\Bigr),
\]
\[
\dot{L}_n=\frac{1}{2(uk+\mathrm{h}^\vee)}
\Bigl(\sum_{m\in\mathbb{Z}}\sum_i
u_i^{(-um)}u^{i(um+un)}\Bigr),
\quad n\neq0,
\]
and $\dot{Z}$ acts by the scalar
\[
c_{uk}=\frac{uk\,\dim(\overline{\mathfrak{g}})}{uk+\mathrm{h}^\vee}.
\]

\begin{remark}
The map $\psi_u$ extends to an algebra homomorphism
\[
\psi_u:\mathcal{U}(\mathfrak{g})\to \mathcal{U}(\mathfrak{g}).
\]
Up to a completion allowing infinite series, we have
$\psi_u(L_n)=\dot{L}_n$ and $\psi_u(Z)=\dot{Z}$.
\end{remark}

\begin{proposition}\label{Ln1u}
For any $n,m \in \mathbb{Z}$ and $x \in \overline{\mathfrak{g}}$, we have
\begin{enumerate}
    \item[1.] $[x^{(um)},\dot{L}_n] = m x^{(um+un)}$.
    \item[2.] $[\dot{L}_m,\dot{L}_n]=(m-n)\dot{L}_{m+n}+\delta_{m+n,0}\frac{m^3-m}{12}\dot{Z}$.
\end{enumerate}
\end{proposition}

The map $\zeta_u: \Vir \rightarrow \Vir$, defined by
\[
\zeta_u(L_n)=u^{-1}L_{un}+\delta_{n,0}\frac{u-u^{-1}}{24}Z,
\qquad
\zeta_u(Z)=uZ,
\]
is a Lie algebra homomorphism. Set $\widetilde{L}_n=\zeta_u(L_n)$ and
$\widetilde{Z}=\zeta_u(Z)$.

\begin{proposition}
For any $n,m \in \mathbb{Z}$ and $x \in \overline{\mathfrak{g}}$, we have
\begin{enumerate}
    \item[1.] $[x^{(um)},\widetilde{L}_n]=m x^{(um+un)}$.
    \item[2.] $[\widetilde{L}_m,\widetilde{L}_n]=(m-n)\widetilde{L}_{m+n}+\delta_{m+n,0}\frac{m^3-m}{12}\widetilde{Z}$.
\end{enumerate}
\end{proposition}

Set $L_n^{[u]}=\widetilde{L}_n-\dot{L}_n$ and $Z^{[u]}=\widetilde{Z}-\dot{Z}$.

\begin{proposition}\label{Lu-untwisted}
We have
\begin{enumerate}
    \item[1.] $[\mathfrak{g}'[u],L_n^{[u]}]=0$.
    \item[2.] $[L^{[u]}_m,L^{[u]}_n]=(m-n)L^{[u]}_{m+n}+\delta_{m+n,0}\frac{m^3-m}{12} Z^{[u]}$.
\end{enumerate}
\end{proposition}
\subsubsection{Twisted case}

In the case $r=2,3$, we have
\[
\mathfrak{g}
=\bigoplus_{j\in\mathbb{Z}} t^j\otimes\overline{\mathfrak{g}}_{\overline{j}}
\oplus \mathbb{C}K \oplus \mathbb{C}d',
\qquad
\mathfrak{g}[u]
=\bigoplus_{j\in\mathbb{Z}} t^{uj}\otimes\overline{\mathfrak{g}}_{\overline{uj}}
\oplus \mathbb{C}K \oplus \mathbb{C}d',
\]
where
$\overline{\mathfrak{g}}=\bigoplus_{\overline{j}\in\mathbb{Z}/r\mathbb{Z}}
\overline{\mathfrak{g}}_{\overline{j}}$
is the $\mathbb{Z}/r\mathbb{Z}$-gradation induced by the diagram automorphism
$\sigma$.
For $\Lambda\in P^k_+$, the $\mathfrak{g}$-module $L(\Lambda)$ is viewed as a
$\mathfrak{g}[u]$-module of level $uk$.
Fix dual bases $\{u_{i,-\overline{j}}\}$ of
$\overline{\mathfrak{g}}_{-\overline{j}}$ and
$\{u^{i,\overline{j}}\}$ of $\overline{\mathfrak{g}}_{\overline{j}}$, and set
$x^{(n)}=t^n\otimes x$. The \black{Virasoro operators} $L_n$ ($n\in\mathbb{Z}$) and $Z$ on $L(\Lambda)$,
given by the Sugawara construction, are
\[
L_0=\frac{1}{2(k+\mathrm{h}^\vee)}\!\left(
\sum_i u_{i,\overline{0}}u^{i,\overline{0}}
+2\!\sum_{m>0}\sum_i u_{i,-\overline{m}}^{(-m)}u^{i,\overline{m}(m)}
+\frac{r-1}{2r}\bigl(\dim\overline{\mathfrak{g}}
-\dim\overline{\mathfrak{g}}_{\overline{0}}\bigr)k\,\mathrm{Id}_V
\right),
\]
\[
L_n=\frac{1}{2(k+\mathrm{h}^\vee)}
\sum_{m\in\mathbb{Z}}\sum_i
u_{i,-\overline{m}}^{(-m)}u^{i,\overline{m}(m+rn)},
\quad n\neq0,
\]
and $Z$ acts by the scalar
\[
 c_k=\frac{rk\,\dim(\overline{\mathfrak{g}})}{k+\mathrm{h}^\vee}.
\]

The \black{Virasoro operators} $\dot{L}_n$ ($n\in\mathbb{Z}$) and $\dot{Z}$ on
$\mathfrak{g}[u]$-modules are defined by
\[
\dot{L}_0=\frac{1}{2(uk+\mathrm{h}^\vee)}\!\left(
\sum_i u_{i,\overline{0}}u^{i,\overline{0}}
+2\!\sum_{m>0}\sum_i
u_{i,-\overline{um}}^{(-um)}u^{i,\overline{um}(um)}
+\frac{r-1}{2r}\bigl(\dim\overline{\mathfrak{g}}
-\dim\overline{\mathfrak{g}}_{\overline{0}}\bigr)uk\,\mathrm{Id}_V
\right),
\]
\[
\dot{L}_n=\frac{1}{2(uk+\mathrm{h}^\vee)}
\sum_{m\in\mathbb{Z}}\sum_i
u_{i,-\overline{um}}^{(-um)}u^{i,\overline{um}(um+urn)},
\quad n\neq0,
\]
and $\dot{Z}$ acts by the scalar
\[
c_{uk}=\frac{urk\,\dim(\overline{\mathfrak{g}})}{uk+\mathrm{h}^\vee}.
\]

\begin{proposition}\label{Lnru}
For any $n,m \in \mathbb{Z}$ and $x \in \overline{\mathfrak{g}}_{\overline{m}}$, we have
\begin{enumerate}
    \item[1.] $[x^{(um)},\dot{L}_n] = \frac{m}{r} x^{(um+urn)}$.
    \item[2.] $[\dot{L}_m,\dot{L}_n]
    =(m-n)\dot{L}_{m+n}
    +\delta_{m+n,0}\frac{m^3-m}{12}\dot{Z}$.
\end{enumerate}
\end{proposition}

The map $\zeta_u: \Vir \rightarrow \Vir$, defined by
\[
\zeta_u(L_n)=u^{-1}L_{un}+\delta_{n,0}\frac{u-u^{-1}}{24}Z,
\qquad
\zeta_u(Z)=uZ,
\]
is a Lie algebra homomorphism.
Set $\widetilde{L}_n=\zeta_u(L_n)$ and $\widetilde{Z}=\zeta_u(Z)$.

\begin{proposition}\label{w0LZut}
For any $n,m \in \mathbb{Z}$ and $x \in \overline{\mathfrak{g}}_{\overline{um}}$, we have
\begin{enumerate}
    \item[1.] $[x^{(um)},\widetilde{L}_n]=m x^{(um+urn)}$.
    \item[2.] $[\widetilde{L}_m,\widetilde{L}_n]
    =(m-n)\widetilde{L}_{m+n}
    +\delta_{m+n,0}\frac{m^3-m}{12}\widetilde{Z}$.
\end{enumerate}
\end{proposition}

Set $L_n^{[u]}=\widetilde{L}_n-r\dot{L}_n$ and
$Z^{[u]}=\widetilde{Z}-r\dot{Z}$.

\begin{proposition}\label{Lu-twisted}
We have
\begin{enumerate}
    \item[1.] $[\mathfrak{g}'[u],L_n^{[u]}]=0$.
    \item[2.] $[L^{[u]}_m,L^{[u]}_n]
    =(m-n)L^{[u]}_{m+n}
    +\delta_{m+n,0}\frac{m^3-m}{12} Z^{[u]}$.
\end{enumerate}
\end{proposition}
\subsubsection{Anti-involution $\omega_0$ on Virasoro operators}

The first conclusions of Propositions \ref{Ln1u} and \ref{Lnru} allow us to define a
$\mathfrak{g}'[u]\rtimes \Vir$-module structure on $L(\Lambda)$. Here $\Vir$ denotes
the Virasoro algebra with Virasoro operators $\dot{L}_n$ ($n\in \mathbb{Z}$) and
$\dot{Z}$. We recall the anti-linear anti-involution $\omega_0$ on
$\mathcal{U}(\mathfrak{g})$ defined in Subsection~\ref{w0}. Up to completion of
$\mathcal{U}(\mathfrak{g})$ to allow series in place of finite sums, we have the
following statement.

\begin{proposition}\label{w0LZu}
We have $\omega_0(\dot{L}_n)=\dot{L}_{-n}$ for all $n \in \mathbb{Z}$, and
$\omega_0(\dot{Z})=\dot{Z}$.
\end{proposition}

\begin{remark}
Proposition~\ref{w0LZu} implies
\begin{equation*}
    \langle \dot{L}_n v, w \rangle_\Lambda
    = \langle v, \omega_0(\dot{L}_n) w \rangle_\Lambda
    = \langle v, \dot{L}_{-n} w \rangle_\Lambda,
\end{equation*}
and
\begin{equation*}
    \langle \dot{Z} v, w \rangle_\Lambda
    = \langle v, \omega_0(\dot{Z}) w \rangle_\Lambda
    = \langle v, \dot{Z} w \rangle_\Lambda.
\end{equation*}
\end{remark}

By Propositions \ref{w0LZ} and \ref{w0LZu}, we obtain the following result.

\begin{proposition}\label{w0LZuu}
We have $\omega_0(L_n^{[u]})=L^{[u]}_{-n}$ and $\omega_0(Z^{[u]})=Z^{[u]}$.
\end{proposition}

\begin{remark}\label{w0LnZuuvw}
Proposition~\ref{w0LZuu} implies
\begin{equation*}
    \langle L_n^{[u]} v, w \rangle_\Lambda
    = \langle v, \omega_0(L^{[u]}_n) w \rangle_\Lambda
    = \langle v, L_{-n}^{[u]} w \rangle_\Lambda,
\end{equation*}
and
\begin{equation*}
    \langle Z^{[u]} v, w \rangle_\Lambda
    = \langle v, \omega_0(Z^{[u]}) w \rangle_\Lambda
    = \langle v, Z^{[u]} w \rangle_\Lambda.
\end{equation*}
\end{remark}

\addcontentsline{toc}{section}{References}
\bibliography{references}{}
\bibliographystyle{alpha}

\noindent Université de Lyon, Université Lyon~1, CNRS, Institut Camille Jordan UMR~5208,  
F-69622 Villeurbanne, France\\
E-mail: \href{mailto:khanh.mathematic@gmail.com}{khanh.mathematic@gmail.com}

\end{document}